\documentclass[
fontsize=10pt,
DIV=classic,
abstract=true,
bibliography=totoc
]{scrartcl}

\addtokomafont{section}{\centering\mathversion{bold}}
\addtokomafont{subsection}{\normalsize\mathversion{bold}}
\addtokomafont{subsubsection}{\mathversion{bold}}

\RedeclareSectionCommands[
runin=false,
afterindent=false,
beforeskip=.5\baselineskip,
afterskip=.5\baselineskip
]
{section}

\RedeclareSectionCommands[
runin=true,
beforeskip=0\baselineskip,
afterskip=-5pt plus -1sp minus 1sp
]
{subsection}

\RedeclareSectionCommands[
runin=true,
beforeskip=0\baselineskip,
afterskip=-5pt plus -1sp minus 1sp
]
{paragraph}

\usepackage[main=english]{babel}

\usepackage{ellipsis}

\usepackage[onehalfspacing]{setspace}

\usepackage[tracking=true]{microtype}		
\DeclareMicrotypeSet*[tracking]{my}{ font = */*/*/sc/* }
\SetTracking{ encoding = *, shape = sc }{ 45 }		

\usepackage[all,defaultlines=2]{nowidow}
\usepackage[T1]{fontenc}
\usepackage[sc]{mathpazo}

\usepackage{bm}

\usepackage{amsmath,amsfonts,amssymb,amsthm}
\usepackage{mathtools,txfonts,mathrsfs,yhmath,extarrows,tikz-cd,babel,etoolbox}

\forcsvlist{\listadd\UpperAlphabet}{A,B,C,D,E,F,G,H,I,J,K,L,M,N,O,P,Q,R,S,T,U,V,W,X,Y,Z}
\forcsvlist{\listadd\LowerAlphabet}{a,b,c,e,d,e,f,h,i,j,k,l,m,o,p,q,r,s,t,u,v,w,x,y,z}
\forcsvlist{\listadd\LowerAlphabetWithoutI}{a,b,c,e,d,e,f,h,j,k,l,m,o,p,q,r,s,t,u,v,w,x,y,z}

\newcommand{\BlackboardBold}[1]{\expandafter\def\csname b#1\endcsname{\mathbb{#1}}}
\newcommand{\Caligraphy}[1]{\expandafter\def\csname c#1\endcsname{\mathcal{#1}}}
\newcommand{\Fraktur}[1]{\expandafter\def\csname f#1\endcsname{\mathfrak{#1}}}
\newcommand{\Script}[1]{\expandafter\def\csname s#1\endcsname{\mathscr{#1}}}
\newcommand{\Upright}[1]{\expandafter\def\csname r#1\endcsname{\mathrm{#1}}}

\forlistloop\BlackboardBold{\UpperAlphabet}
\forlistloop\Caligraphy{\UpperAlphabet}
\forlistloop\Fraktur{\UpperAlphabet}
\forlistloop\Script{\UpperAlphabet}
\forlistloop\Upright{\UpperAlphabet}

\forlistloop\Fraktur{\LowerAlphabetWithoutI}
\forlistloop\Upright{\LowerAlphabet}

\usepackage[shortlabels]{enumitem}

\newlist{enumroman}{enumerate}{3}
\setlist[enumroman]{label=(\roman*),ref =(\roman*)}
\newlist{enumalph}{enumerate}{3}
\setlist[enumalph]{label=(\alph*),ref=(\alph*)}

\newcommand\restrict[2]{\left.#1\right|_{#2}}

\newcommand\longto{\ensuremath{\longrightarrow}}
\newcommand\longfrom{\ensuremath{\longleftarrow}}

\newcommand\op{\ensuremath{\mathrm{op}}}
\newcommand\wh{\ensuremath{\mathrm{wh}}}

\newcommand\blank{{-}}

\newcommand\Dcap{\wideparen{\cD}}

\AtBeginDocument{
\let\bemph\textbf
\let\del\partial
}

\usepackage{appendix}

\usepackage{hyperref}
\addto\extrasenglish{%
}

\newcommand{\autorefpart}[2]{\hyperref[#2]{\autoref*{#1}\,\ref*{#2}}}
\newcommand{\autorefparts}[3]{\hyperref[#3]{\autoref*{#1}\,\ref*{#2}\,\ref*{#3}}}
\newcommand{\equationref}[1]{\hyperref[#1]{(\ref*{#1})}}

\hypersetup{colorlinks=true,linkcolor=teal,citecolor=violet}

\makeatletter
\patchcmd{\hyper@makecurrent}{%
    \ifx\Hy@param\Hy@chapterstring
        \let\Hy@param\Hy@chapapp
    \fi
}{%
    \iftoggle{inappendix}{%
        \@checkappendixparam{chapter}%
        \@checkappendixparam{section}%
        \@checkappendixparam{subsection}%
        \@checkappendixparam{subsubsection}%
        \@checkappendixparam{paragraph}%
        \@checkappendixparam{subparagraph}%
    }{}%
}{}{\errmessage{failed to patch}}

\newcommand*{\@checkappendixparam}[1]{%
    \def\@checkappendixparamtmp{#1}%
    \ifx\Hy@param\@checkappendixparamtmp
        \let\Hy@param\Hy@appendixstring
    \fi
}
\makeatletter

\newtoggle{inappendix}
\togglefalse{inappendix}

\apptocmd{\appendix}{\toggletrue{inappendix}}{}{\errmessage{failed to patch}}
\apptocmd{\subappendices}{\toggletrue{inappendix}}{}{\errmessage{failed to patch}}

\makeatletter
\newcommand{\leqnomode}{\tagsleft@true\let\veqno\@@leqno}
\makeatother

\usepackage{thmtools}

\declaretheoremstyle[
spaceabove=3pt, spacebelow=3pt,
postheadspace=1ex,
bodyfont=\itshape,
numberwithin=section
]{myplain}

\declaretheoremstyle[
spaceabove=3pt, spacebelow=3pt,
postheadspace=1ex,
numberwithin=section
]{mydefinition}

\declaretheoremstyle[
spaceabove=3pt, spacebelow=3pt,
postheadspace=1ex,
headfont=\itshape,
numberwithin=section
]{myremark}

\declaretheoremstyle[
spaceabove=3pt, spacebelow=3pt,
postheadspace=1ex,
headfont=\itshape,
qed=\qedsymbol
]{myproof}

\declaretheorem[style=myplain,numbered=yes,refname={Theorem,Theorems}]{Theorem}

\declaretheorem[style=myplain,numbered=yes,name=Theorem,refname={Theorem,Theorems}]{MainTheorem}

\declaretheorem[sibling=Theorem,numbered=yes,style=myplain,refname={Proposition,Propositions}]{Proposition}

\declaretheorem[sibling=Theorem,numbered=yes,style=myplain,refname={Lemma,Lemmata}]{Lemma}

\declaretheorem[sibling=Theorem,numbered=yes,style=myplain,refname={Corollary,Corollaries}]{Corollary}

\declaretheorem[sibling=Theorem,numbered=yes,style=mydefinition,refname={Definition,Definitions}]{Definition}

\declaretheorem[sibling=Theorem,numbered=yes,style=myremark,refname={Remark,Remarks}]{Remark}

\declaretheorem[unnumbered,style=myproof]{Proof}

\newcommand{\DeclareMyOperator}[1]{\expandafter\DeclareMathOperator\csname#1\endcsname{#1}}
\forcsvlist{\DeclareMyOperator}{Ann,Ch,coker,Der,Ext,gr,Gr,id,Proj,Sp,Spec,Supp,Sym,HH}

\DeclareMathOperator{\iHom}{\mathcal{H}\kern-.1em\mathit{om}}

\newcommand{\DeclareMyOperatorStar}[1]{\expandafter\DeclareMathOperator*\csname#1\endcsname{#1}}
\forcsvlist{\DeclareMyOperatorStar}{}

\newcommand{\DeclareMyCategory}[1]{\expandafter\def\csname#1\endcsname{{\mathsf{#1}}}}
\forcsvlist{\DeclareMyCategory}{Coh,Hol}

\usepackage{titling}
\predate{}
\postdate{}

\title{\sffamily\bfseries On \(\bm\Dcap\)-Modules of Finite Length on Rigid Analytic Spaces}
\author{Julian Reichardt}
\date{}

\begin{document}

\maketitle

\begin{abstract}\noindent
We show that for quasi-compact smooth rigid analytic spaces, the extension functor sends holonomic \(\cD\)-modules to coadmissible \(\Dcap\)-modules which are of finite length as weakly holonomic \(\Dcap\)-modules. Using this, we show that the meromorphic connections considered by Bitoun--Bode and the local cohomology groups considered by Ardakov--Bode--Wadsley are of finite length as weakly holonomic \(\Dcap\)-modules for quasi-compact smooth rigid analytic spaces. As a central tool, we introduce and study Hilbert polynomials for finitely generated modules over completed Weyl algebras.
\end{abstract}

\section{Introduction}

The theory of algebraic \(\cD\)-modules on smooth complex varieties is a classical topic, with numerous important applications to algebraic geometry and geometric representation theory. In \cite{AW2019}, K. Ardakov and S. Wadsley initiated the study of an analogous framework of so-called \(\Dcap\)-modules for smooth rigid analytic varieties over a nonarchimedean base field \(K\) of mixed characteristic \((0,p)\). The starting point is the construction of a suitable completion of the sheaf of algebraic differential operators to account for the ubiquity of algebro-topological structures one encounters throughout rigid geometry. 

Although the theory of \(\Dcap\)-modules has been extensively developed further by various contributors, certain foundational aspects remain difficult to approach. Of particular notice is the development of a robust notion of holonomicity for \(\Dcap\)-modules. Recall that holonomic \(\cD\)-modules on smooth complex varieties can be characterised as those coherent \(\cD\)-modules whose associated characteristic variety is of minimal dimension (see \cite[Section~3]{HTT2008}). In contrast to arbitrary coherent \(\cD\)-modules, holonomic \(\cD\)-modules enjoy stability under the standard operations for \(\cD\)-modules and are always of finite length as \(\cD\)-modules. Note, however, that not every \(\cD\)-module of finite length is holonomic. 

In the absence of a suitable replacement for the characteristic variety, Ardakov--Wadsley, jointly with A.\ Bode, introduced weakly holonomic \(\Dcap\)-modules in \cite{ABW2021}. These are defined in terms of a purely homological characterisation of holonomic \(\cD\)-modules, which can be suitably translated to \(\Dcap\)-modules. However, there are pathological examples of weakly holonomic \(\Dcap\)-modules that are of infinite length or have infinite dimensional fibres (see \cite[Section~8.1]{ABW2021}). As an alternative, Bode proposed in \cite{Bode2025Holonomic} a definition of a (derived) category of holonomic \(\Dcap\)-modules, taking stability under the standard operations essentially as a definition. This is inspired by the work of D. Caro \cite{Caro2004,Caro2009} on overholonomic arithmetic \(\sD^\dagger\)-modules. This produces a functorially well-behaved category of holonomic \(\Dcap\)-modules containing integrable connections. However, due to its abstract nature, both the precise relation to the previously considered weakly holonomic \(\Dcap\)-modules as well as questions of finite length are a priori unclear.

In this article, we develop a framework providing a criterion for weakly holonomic \(\Dcap\)-modules on the closed unit polydisc to be of finite length. More precisely, given a coadmissible presentation \(M=\varprojlim_n M_n\) of a weakly holonomic \(\Dcap\)-module \(M\) on the closed unit polydisc, we give a criterion for \(M\) to be of finite length by considering uniform bounds on the lengths of the \(M_n\) (see \autoref{Proposition:Finite_Length_Criterion}).

As an application, our first main result is concerned with the extension functor from coherent \(\cD\)-modules to coadmissible \(\Dcap\)-modules. This functor has been considered before in \cite{BB2021} and \cite{ABW2021} and is a valuable tool for studying \(\Dcap\)-modules. It is known that the extension functor sends holonomic \(\cD\)-modules to weakly holonomic \(\Dcap\)-modules. We obtain the following strengthening.

\begin{MainTheorem}[\autoref{Theorem:Finite_Length_Affinoid}]\label{MainTheorem:Extension_Finite_Length}
Let \(X=\Sp A\) be a smooth affinoid \(K\)-space and let \(\cM\) be a holonomic \(\cD_X\)-module. Then the extended \(\Dcap_X\)-module \(E_X(\cM)=\Dcap_X\otimes_{\cD_X}\cM\) is of finite length in the category of weakly holonomic \(\Dcap_X\)-modules. 
\end{MainTheorem}

Note that \autoref{MainTheorem:Extension_Finite_Length} immediately extends to quasi-compact smooth rigid analytic \(K\)-spaces by a simple covering argument (see \autoref{Corollary:Finite_Length_Quasi_Compact}). In some sense, \autoref{MainTheorem:Extension_Finite_Length} is a purely algebraic statement. However, the extension functor can be used to study questions of coadmissibility and weak holonomicity of \(\Dcap\)-modules, as employed before in \cite{BB2021} and \cite{ABW2021}. This way, we show that various weakly holonomic \(\Dcap\)-modules considered before are of finite length. 

\begin{MainTheorem}[\autoref{Theorem:ABW_Open_Extension}]\label{MainTheorem:Geometric_Finite_Length:Pushforward}
Let \(X\) be a quasi-compact smooth rigid analytic \(K\)-space. Let \(j\colon U\to X\) be the embedding of a Zariski-open subspace. Let \(\cM\) be an integrable connection on \(X\). Then \(\rR^qj_\ast\bigl(\restrict{\cM}{U}\bigr)\) is of finite length in the category of weakly holonomic \(\Dcap_X\)-modules for all \(q\ge0\).
\end{MainTheorem}

\begin{MainTheorem}[\autoref{Theorem:ABW_Closed_Extension}]\label{MainTheorem:Geometric_Finite_Length:Cohomology}
Let \(X\) be a quasi-compact smooth rigid analytic \(K\)-space. Let \(Z\) be a closed analytic subset and let \(\cM\) be an integrable connection on \(X\). Then \(\underline\cH_Z^q(\cM)\) is of finite length in the category of weakly holonomic \(\Dcap_X\)-modules for all \(q\ge0\).
\end{MainTheorem}

In the case where \(X\) is a smooth affinoid \(K\)-space with free tangent sheaf and the complement \(Z=X\setminus U\) is defined by a single equation, we obtain a more refined version of \autoref{MainTheorem:Geometric_Finite_Length:Pushforward} for meromorphic connections on \(X\) with singularities along \(Z\). This is a strengthening of \cite[Theorem~4.4]{BB2021}. As the statement is quite technical, we postpone it to the body of the text (see \autoref{Theorem:BB_Extension}). In particular, this provides explicit examples of \(\Dcap\)-modules of finite length with non-trivial monodromy. Moreover, this statement can be used to obtain examples of non-holonomic weakly holonomic \(\Dcap\)-modules of finite length (see \autoref{Corollary:Non_Holonomic_Finite_Length}). As in the classical setting, this shows that finite length is not sufficient for holonomicity.

\paragraph{Structure of the paper} \autoref{Section:Stabilisation} and \autoref{Section:Completed_Weyl_Algebras} constitute the technical foundation of this paper. 

We begin in \autoref{Section:Stabilisation} by studying deformations of modules over deformable \(R\)-algebras and doubly filtered \(K\)-algebras, as considered in \cite{AW2013,ABW2021}. Our main results are the stabilisation results \autoref{Proposition:Stabilisation_Deformable_Algebra} and \autoref{Proposition:Stabilisation_Doubly_Filtered_Algebra}. 

In \autoref{Section:Completed_Weyl_Algebras} we study a family of completions of the Weyl algebra over \(K\) as doubly filtered \(K\)-algebras. We introduce double filtrations based on the Bernstein filtration, which allow the definition of a Hilbert polynomial for finitely generated modules (see \autoref{Definition:Hilbert_Polynomial_Completed_Weyl_Algebras}). As an application, we show that holonomic modules over these completed Weyl algebras are of finite length, in analogy to holonomic modules over the Weyl algebra over \(K\) (see \autoref{Corollary:Multiplicity_Length_Bound}).

Moving towards geometry, we use \autoref{Section:Differential_Operators} to survey the theory of algebraic \(\cD\)-modules and analytic \(\Dcap\)-modules on smooth rigid analytic \(K\)-spaces. Here we approach both theories uniformly using the language of Lie--Rinehart algebras, the algebraic incarnations of Lie algebroids. Of particular importance is \autoref{Subsection:Dcap_Unit_Polydisc}, where we make the theory of \(\Dcap\)-modules on the closed unit polydisc explicit in terms of the completed Weyl algebras we studied in \autoref{Section:Completed_Weyl_Algebras}. This also prepares the applications of our stabilisation results from \autoref{Section:Stabilisation}.

We prove our main results in \autoref{Section:Extension} and \autoref{Section:Finite_Length}.

\autoref{Section:Extension} is devoted to the proof of \autoref{MainTheorem:Extension_Finite_Length} (see \autoref{Theorem:Finite_Length_Affinoid}). We first recall the definition of the extension functor and establish some compatibility results. Using Kashiwara's equivalence, we then reduce the general affinoid case to a simplified case on the closed unit polydisc. Here we can apply the machinery developed in \autoref{Section:Stabilisation} and \autoref{Section:Completed_Weyl_Algebras}.

Finally, in \autoref{Section:Finite_Length} we prove \autoref{MainTheorem:Geometric_Finite_Length:Pushforward} and \autoref{MainTheorem:Geometric_Finite_Length:Cohomology} (see \autoref{Theorem:ABW_Open_Extension} and \autoref{Theorem:ABW_Closed_Extension}). As outlined before, this is done by reducing to the previously established theorems for extensions of holonomic \(\Dcap\)-modules. We also obtain the more refined version \autoref{Theorem:BB_Extension} of \autoref{Theorem:ABW_Open_Extension}, based on \cite{BB2021}, and show that there are non-holonomic weakly holonomic \(\Dcap\)-modules of finite length in \autoref{Corollary:Non_Holonomic_Finite_Length}.

This paper has two appendices. In \autoref{Appendix:Hilbert_Polynomials} we recall the construction and fundamental properties of Hilbert polynomials for polynomial rings over general base fields. In \autoref{Appendix:D_Modules_Rigid_Spaces} we collect some generalities on coordinate systems for smooth affinoid \(K\)-spaces and go over some subtleties which arise when working with algebraic \(\cD\)-modules on smooth rigid analytic \(K\)-spaces.

\paragraph{Acknowledgments} This work is part of the author's PhD thesis, supervised by Tobias Schmidt. I would like to thank him for gently guiding me into this topic. I would also like to thank Andreas Bode, for many enlightening discussions, and Feliks R{\k{a}}czka, for answering questions about \cite{Raczka2024Thesis}. Moreover, I wish to express my gratitude to Tobias Schmidt and Andreas Bode for carefully reading a first draft of this article. This work was carried out as part of the research training group \emph{GRK 2240: Algebro-Geometric Methods in Algebra, Arithmetic and Topology} funded by the \emph{Deutsche Forschungsgemeinschaft}.

\paragraph{Notations and Conventions} We fix a complete discrete valuation field \(K\) of mixed characteristic \((0,p)\) with ring of integers \(R\) and residue field \(k\). We also fix a uniformiser \(\pi\). For an \(R\)-module \(M\), we denote by \(\widehat{M}\) its \(\pi\)-adic completion and by \(M_K\) its extension of scalars from \(R\) to \(K\). We also abbreviate \(\widehat{M}\otimes_R K\) to \(\widehat{M_K}\). Rings are assumed to be associative and unital. By noetherian we mean left and right noetherian and similar for other ring-theoretic adjectives. Modules will be left modules if not specified further. We work with rigid analytic \(K\)-spaces in the sense of Tate (see \cite{BGR1984}). A rigid analytic \(K\)-space \(X\) will always assumed to be quasi-separated, quasi-paracompact and equidimensional.

\section{Stabilisation of Associated Graded Modules}\label{Section:Stabilisation}

This section is concerned with the formalisms of deformable \(R\)-algebras and doubly filtered \(K\)-algebras, as introduced in \cite[Section~3]{AW2013}. Our primary goal is to adapt the \emph{graded stabilisation} results \cite[Lemma~3.5]{AW2013} for deformable \(R\)-algebras to deformations of finitely generated modules, as discussed in \cite[Section~3]{ABW2021}. This also implies analogous stabilisation results for the families of doubly filtered \(K\)-algebras associated to a deformable \(R\)-algebra.

\subsection{Graded stabilisation for deformable \texorpdfstring{\(R\)}{R}-algebras} Let \(U\) be a positively filtered \(R\)-algebra such that \(F_0U\) is an \(R\)-subalgebra. Recall from \cite[Section~3.5]{AW2013}, that we call \(U\) a \bemph{deformable}, if \(\gr U\) is a flat \(R\)-module. The \bemph{\(\bm n\)-th deformation of \(\bm U\)} is the \(R\)-submodule
\[
\bm{U_n}
\coloneqq\sum_{i\ge0}\pi^{in}F_i U\,.
\]     
Note that \(U_n\) is an \(R\)-subalgebra of \(U\). As \(\gr U\) is flat over \(R\), this subspace filtration on \(U_n\) takes the form
\[
F_i U_n
=\sum_{j=0}^i\pi^{jn}F_j U\,.
\]
We have the following stabilisation result on the associated graded \(R\)-algebras of the deformations of a deformable \(R\)-algebra.

\begin{Proposition}[{{\cite[Lemma~3.5]{AW2013}}}]\label{Proposition:Graded_Isomorphism_Deformed_Algebra}
Let \(U\) be a deformable \(R\)-algebra. Then \(U_n\) is a deformable \(R\)-algebra and the natural morphism
\[
\xi_n\colon\gr U\longto\gr U_n
,\qquad
u+F_{i-1}U\longmapsto \pi^{in}u+F_{i-1}U_n
\]
is a graded isomorphism of \(R\)-algebras for all \(n\ge0\).
\end{Proposition}

If \(U\) is a deformable \(R\)-algebra, then there is an isomorphism \((U_n)_m\cong U_{n+m}\) of \(R\)-algebras by \cite[Lemma~6.4(b)]{AW2019}. This yields the following consequence of \autoref{Proposition:Graded_Isomorphism_Deformed_Algebra}.

\begin{Corollary}\label{Corollary:Iterated_Deformations}
Let \(U\) be a deformable \(R\)-algebra and let \(n_0\ge0\). The natural morphism
\[
\xi_{n_0,n}\colon\gr U_{n_0}\longto\gr U_n
,\qquad
u+F_{i-1}U_{n_0}\longmapsto \pi^{i(n-n_0)}u+F_{i-1}U_n
\]
is a graded isomorphism of \(R\)-algebras satisfying \(\xi_n=\xi_{n_0,n}\xi_{n_0}\) for all \(n\ge n_0\).
\end{Corollary}

\begin{Corollary}\label{Corollary:Deformable_Flat}
Let \(U\) be a deformable \(R\)-algebra. Then \(U_n\) is a flat \(R\)-module for all \(n\ge0\).
\end{Corollary}

\begin{Proof}
As \(R\) is a discrete valuation ring by assumption, an \(R\)-module \(M\) is flat if and only if the \(\pi\)-multiplication map \(\pi_M\colon M\to M,\,m\mapsto\pi m\) is injective. Since \(U_n\) is deformable for each \(n\ge0\) by \autoref{Proposition:Graded_Isomorphism_Deformed_Algebra}, it suffices to check that the deformable \(R\)-algebra \(U\) is flat as an \(R\)-module. As \(F_0U\) is an \(R\)-subalgebra, \(\pi_U\) is filtered and induces \(\pi_{\gr U}\). As \(\gr U\) is flat \(R\)-module, \(\pi_{\gr U}\) is injective. Thus, \(\pi_U\) is injective by \cite[Corollary~I.4.2.5(1)]{HvO1996}.
\end{Proof}

Let \(U\) be a deformable \(R\)-algebra and let \(M\) be a finitely generated \(U\)-module. Fix a good filtration \(F_\bullet M\) on \(M\) (e.g., for \(m_1,\dots,m_r\) a finite generating set of \(M\), consider the filtration \(F_iM=\sum_{j=1}^r(F_i U)m_j\); cf.\ \cite[Remark~I.5.2]{HvO1996}). The \bemph{\(\bm n\)-th deformation of \(\bm M\)} is the finitely generated \(U_n\)-submodule
\[
\bm{M_n}
\coloneqq \sum_{i\ge0}\pi^{in}F_iM\,.
\]
On \(M_n\) we have the following natural filtration
\[
F_iM_n
=\sum_{j=0}^i\pi^{jn}F_j M\,.
\]
If \(\gr M\) is flat over \(R\) (equivalently, \(\pi\)-torsionfree), then this filtration agrees with the subspace filtration on \(M_n\). We have the following partial analogue of \autoref{Proposition:Graded_Isomorphism_Deformed_Algebra} for the deformations of finitely generated \(U\)-modules, whose proof is immediate from the definitions.

\begin{Lemma}[{{\cite[Lemma~3.1]{ABW2021}}}]\label{Lemma:Graded_Morphism_Deformed_Module}
Let \(U\) be a deformable \(R\)-algebra and let \(M\) be a \(\pi\)-torsionfree finitely generated \(U\)-module. The natural morphism 
\[
\mu_n\colon\gr M\longto\gr M_n
,\qquad
m+F_{i-1}M\longmapsto \pi^{in}m+F_{i-1}M_n
\]
is a graded \(R\)-linear morphism such that
\begin{enumroman}
\item\label{Lemma:Graded_Morphism_Deformed_Module:Linearity} \(\mu_n(um)=\xi_n(u)\mu_n(m)\) for all \(u\in U\) and \(m\in\gr M\),
\item\label{Lemma:Graded_Morphism_Deformed_Module:Surjectivity} \(\mu_n\) is surjective, and
\item\label{Lemma:Graded_Morphism_Deformed_Module:Kernel} \(\ker\bigl(\restrict{\mu_n}{\gr_i M}\bigr)=\bigl(\gr_i M\bigr)[\pi^{in}]\)
\end{enumroman}
for all \(n\ge0\).
\end{Lemma}

Note that in contrast to the \(\xi_n\) of \autoref{Proposition:Graded_Isomorphism_Deformed_Algebra}, the \(\mu_n\) of \autoref{Lemma:Graded_Morphism_Deformed_Module} are not necessarily injective. However, we have an analogue of \autoref{Corollary:Deformable_Flat}. For a deformable \(R\)-algebra \(U\), we consider the following assumptions:
\begin{equation}\label{Assumption}\leqnomode\tag{\(A\)}
\parbox{\dimexpr\linewidth-4em}{\itshape \(\gr U\) is commutative noetherian \(R\)-algebra and \(F_0 U\) is \(\pi\)-adically separated}
\end{equation}

\begin{Lemma}\label{Lemma:Torsion_Deformed_Module}
Let \(U\) be a deformable \(R\)-algebra satisfying \eqref{Assumption} and let \(M\) be a \(\pi\)-torsionfree finitely generated \(U\)-module. For any good filtration on \(M\), the \(\gr U\)-module \(\gr M\) has bounded \(\pi\)-torsion and there exists an integer \(n_0\ge0\) such that \(\gr M_n\) is \(\pi\)-torsionfree for all \(n\ge n_0\).
\end{Lemma}

\begin{Proof}
Under the given assumptions, \cite[Corollary~3.1]{ABW2021} shows that \(\gr M_n\) is \(\pi\)-torsionfree if \(\pi^n\) annihilates the \(\pi\)-torsion of \(\gr M\). Consider the \(\pi\)-torsion submodule \(\bigl(\gr M\bigr)[\pi^\infty]\) of \(\gr M\). Since \(\gr U\) is noetherian and \(\gr M\) is finitely generated, it follows that \(\bigl(\gr M\bigr)[\pi^\infty]\) is finitely generated as well. In particular, there is \(n_0\) such that \(\pi^{n_0}\bigl(\gr M\bigr)[\pi^\infty]=0\), meaning \(\pi^{n_0}\) annihilates the \(\pi\)-torsion of \(\gr M\). Thus, it follows from \cite[Corollary~3.1]{ABW2021} that \(\gr M_n\) is \(\pi\)-torsionfree for all \(n\ge n_0\).
\end{Proof}

In combining \autoref{Lemma:Graded_Morphism_Deformed_Module} and \autoref{Lemma:Torsion_Deformed_Module}, we obtain an direct analogue of  \autoref{Corollary:Iterated_Deformations} for deformations of \(\pi\)-torsionfree finitely generated \(U\)-modules. 

\begin{Proposition}\label{Proposition:Stabilisation_Deformable_Algebra}
Let \(U\) be a deformable \(R\)-algebra satisfying \eqref{Assumption} and let \(M\) be a \(\pi\)-torsionfree finitely generated \(U\)-module. For any good filtration on \(M\), there exists an integer \(n_0\ge0\) so that the natural morphism
\[
\mu_{n_0,n}\colon\gr M_{n_0}\longto\gr M_n
,\qquad
m+F_{i-1}M_{n_0}\longmapsto\pi^{i(n-n_0)}m+F_{i-1}M_n
\]
is a graded \(R\)-linear isomorphism satisfying \(\mu_n=\mu_{n_0,n}\mu_{n_0}\) for all \(n\ge n_0\).

Moreover, this isomorphism identifies the \(\gr U_{n_0}\)-module structure on \(\gr M_{n_0}\) with the \(\gr U_n\)-module structure on \(\gr M_n\) along the isomorphism \(\xi_{n_0,n}\colon\gr U_{n_0}\to\gr U_n\) of \autoref{Corollary:Iterated_Deformations}.
\end{Proposition}

\begin{Proof}
Denote the \(\pi\)-torsion submodule of \(\gr M\) by \(\bigl(\gr M\bigr)[\pi^\infty]\). By \autoref{Lemma:Torsion_Deformed_Module}, there is an integer \(n_0\ge0\) such that \(\pi^{n_0}\bigl(\gr M\bigr)[\pi^\infty]=0\) and such that \(\gr M_n\) is \(\pi\)-torsionfree for \(n\ge n_0\).

Now \autorefpart{Lemma:Graded_Morphism_Deformed_Module}{Lemma:Graded_Morphism_Deformed_Module:Linearity} and \autorefpart{Lemma:Graded_Morphism_Deformed_Module}{Lemma:Graded_Morphism_Deformed_Module:Surjectivity} show that there is a surjection \(\mu_n\colon\gr M\to\gr M_n\) for all \(n\ge0\). By construction, the natural morphism induces a factorisation \(\mu_n=\mu_{n_0,n}\mu_{n_0}\) as graded \(R\)-linear morphisms. In particular, \(\mu_{n_0,n}\) is surjective. Moreover, note that \(\mu_{n_0,n}(um)=\xi_{n_0,n}(u)\mu_{n_0,n}(m)\) for all \(u\in\gr U_{n_0}\), for all \(m\in\gr M_{n_0}\) and \(\xi_{n_0,n}\) the natural isomorphism of \autoref{Corollary:Iterated_Deformations}. Thus, it remains to show that \(\mu_{n_0,n}\) is injective to conclude. Restricting the factorisation \(\mu_n=\mu_{n_0,n}\mu_{n_0}\) to the summand \(\gr_i M\) induces a factorisation
\[
\restrict{\mu_n}{\gr_i M}
=\restrict{\mu_{n_0,n}}{\gr_i M_{n_0}}\restrict{\mu_{n_0}}{\gr_i M}\,.
\]
Since \(\mu_{n_0}\) is a graded surjection by \autorefpart{Lemma:Graded_Morphism_Deformed_Module}{Lemma:Graded_Morphism_Deformed_Module:Surjectivity}, the morphism \(\restrict{\mu_{n_0}}{\gr_i M}\) is surjective. Hence the kernel-cokernel sequence of the composition \(\mu_n=\mu_{n_0,n}\mu_{n_0}\) reduces to a short exact sequence
\begin{equation*}\label{Diagram:Kernel_Sequence}
\begin{tikzcd}
0 \ar[r] & \ker\bigl(\restrict{\mu_{n_0}}{\gr_i M}\bigr) \ar[r] & \ker\bigl(\restrict{\mu_n}{\gr_i M}\bigr) \ar[r] & \ker\biggl(\restrict{\mu_{n,n_0}}{\gr_i M_{n_0}}\biggr) \ar[r] & 0\,.\tag{\(\ast\)}
\end{tikzcd}
\end{equation*}
It remains to show that \(\ker\bigl(\restrict{\mu_{n,n_0}}{\gr_i M_{n_0}}\bigr)=0\) for all \(n\ge n_0\) and all \(i\ge0\). By \equationref{Diagram:Kernel_Sequence}, it suffices to show \(\ker\bigl(\restrict{\mu_{n_0}}{\gr_i M}\bigr)=\ker\bigl(\restrict{\mu_n}{\gr_i M}\bigr)\) for all \(n\ge n_0\) and all \(i\ge0\).

As we have chosen \(n_0\) such that \(\pi^{n_0}\bigl(\gr M\bigr)[\pi^\infty]=0\), it follows that 
\[
\bigl(\gr M\bigr)[\pi^n]=\bigl(\gr M\bigr)[\pi^\infty]=\bigl(\gr M\bigr)[\pi^{n_0}]
\]
for all \(n\ge n_0\). Therefore we find that
\[
\bigl(\gr_i M\bigr)[\pi^{n_0i}]
=\bigl(\gr_i M\bigr)\cap\bigl(\gr M\bigr)[\pi^{n_0}]
=\bigl(\gr_i M\bigr)\cap\bigl(\gr M\bigr)[\pi^n]
=\bigl(\gr_i M\bigr)[\pi^{ni}]
\]
for all \(n\ge n_0\) and \(i\ge1\). Hence
\[
\ker\bigl(\restrict{\mu_{n_0}}{\gr_i M}\bigr)
=\bigl(\gr_i M\bigr)[\pi^{n_0i}]
=\bigl(\gr_i M\bigr)[\pi^{ni}]
=\ker\bigl(\restrict{\mu_n}{\gr_i M}\bigr)
\]
by \autorefpart{Lemma:Graded_Morphism_Deformed_Module}{Lemma:Graded_Morphism_Deformed_Module:Kernel} for all \(n\ge n_0\) and \(i\ge1\). Finally, as \(\ker\bigl(\restrict{\mu_n}{\gr_0 M}\bigr)=0\) independently of \(n\), we conclude that \(\ker\bigl(\restrict{\mu_{n_0}}{\gr_i M}\bigr)=\ker\bigl(\restrict{\mu_n}{\gr_i M}\bigr)\) for all \(n\ge n_0\) and all \(i\ge0\). This concludes the proof.
\end{Proof}

\begin{Remark} 
The results of this subsection can be extended to non-discretely valued \(R\). For this, \eqref{Assumption} is replaced by the assumption that \(\gr U\) is commutative and almost noetherian, \(F_0R\) is a commutative \(R\)-algebra topologically of finite type and \(F_i U\) is finitely generated over \(F_0 U\) for \(i\ge0\). Then one can invoke \cite[Lemma~2.11]{Bode2026} instead of \autoref{Lemma:Torsion_Deformed_Module} to carry out the proofs.
\end{Remark}

\subsection{Graded stabilisation for doubly filtered \texorpdfstring{\(K\)}{K}-algebras}\label{Subsection:Stabilisation_Doubly_Filtered_Algebras} Following \cite[Section~2.7]{AW2013}, we say that an \(R\)-submodule \(L\) of a \(K\)-vector space \(V\) is an \bemph{\(\bm R\)-lattice} if \(V=K\cdot L\) and \(\bigcap_{i=0}^\infty\pi^i L=0\). The \(k\)-vector space \(\gr_0 V\coloneqq L/\pi L\) is called the \bemph{slice of \(V\)}. Let \(A\) be a \(K\)-algebra. Recall from \cite[Definition~3.1]{AW2013}, that we say \(A\) is a \bemph{(complete) doubly filtered \(\bm K\)-algebra} if it has an \(R\)-subalgebra \(F_0 A\) which is a (\(\pi\)-adically complete) \(R\)-lattice and if the slice \(\gr_0 A\) of \(A\) admits a (complete) filtration \(F_\bullet\gr_0 A\). For any (complete) doubly filtered \(K\)-algebra \(A\), denote by \(\bm{\Gr(A)}\coloneqq\gr\bigl(\gr_0 A\bigr)\) the \bemph{associated graded \(\bm k\)-algebra}.

Any deformable \(R\)-algebra gives rise to a family of complete doubly filtered \(K\)-algebras.

\begin{Proposition}[{{\cite[Lemma~3.7,~Corollary~3.7]{AW2013}}}]\label{Proposition:Graded_Isomorphism_Doubly_Filtered_Algebra}
Let \(U\) be a deformable \(R\)-algebra. Then \(\widehat{U_{n,K}}\) is a complete doubly filtered \(K\)-algebra with slice \(\gr_0\widehat{U_{n,K}}\cong U_n/\pi U_n\) for all \(n\ge0\). Moreover, the natural morphisms
\[
\Gr\bigl(\widehat{U_{n,K}}\bigr)
\longto\gr\bigl(\widehat{U_n}/\pi\widehat{U_n}\bigr)
\cong\gr\bigl(U_n/\pi U_n\bigr)
\longfrom\gr U_n/\pi\gr U_n
\]
and
\[
\Gr\bigl(\widehat{U_K}\bigr)
\longto\Gr\bigl(\widehat{U_{n,K}}\bigr)
\]
are graded isomorphisms of \(k\)-algebras for all \(n\ge0\).
\end{Proposition}

\begin{Remark}\label{Remark:Iterated_Double_Filtrations}
As in \autoref{Corollary:Iterated_Deformations}, for any deformable \(R\)-algebra \(U\), we obtain natural isomorphisms \(\Gr\bigl(\widehat{U_{n_0,K}}\bigr)\cong\Gr\bigl(\widehat{U_{n,K}}\bigr)\) for all \(n\ge n_0\).
\end{Remark}

We are looking for an analogue of \autoref{Proposition:Graded_Isomorphism_Doubly_Filtered_Algebra} in the spirit of \autoref{Proposition:Stabilisation_Deformable_Algebra}. Let \(A\) be a doubly filtered \(K\)-algebra and let \(M\) be an \(A\)-module. Recall from \cite[Section~3.2]{AW2013}, that a \bemph{double filtration} on \(M\) consists of an \(R\)-lattice \(F_0M\), which is an \(F_0A\)-submodule, and a filtration \(F_\bullet\gr_0M\) on the slice of \(M\) compatible with the filtration \(F_\bullet\gr_0A\) on the slice of \(A\). A double filtration on \(M\) is called \bemph{good}, if the filtration on \(\gr_0 M\) is separated and \(\Gr(M)\) is a finitely generated \(\Gr(A)\)-module. We call \(\bm{\Gr(M)}\coloneqq\gr\bigl(\gr_0M\bigr)\) the \bemph{associated graded \(\bm{\Gr(A)}\)-module}. 

We need the following result from \cite{AW2013} on the existence of good double filtrations for suitable doubly filtered \(K\)-algebras. As we will exploit a subtle point in the construction of such good double filtrations, we repeat the relevant parts of the proof here.

\begin{Lemma}[{{\cite[Proposition~3.2]{AW2013}}}]\label{Lemma:Existence_Double_Filtrations}
Let \(A\) be a complete doubly filtered \(K\)-algebra and assume that \(\Gr(A)\) is noetherian. 
\begin{enumroman}
\item\label{Lemma:Existence_Double_Filtrations:Lattice} Every finitely generated \(\pi\)-torsionfree \(F_0A\)-module \(L\) is an \(R\)-lattice in \(L_K\).
\item\label{Lemma:Existence_Double_Filtrations:Filtration} Every finitely generated \(A\)-module \(M\) has at least one good double filtration.
\end{enumroman}
\end{Lemma}

\begin{Proof}
For the proof of \ref{Lemma:Existence_Double_Filtrations:Lattice} we refer to \cite{AW2013}. As \(\Gr(A)\) is noetherian by assumption, \(F_0A\) is noetherian by \cite[Lemma~3.2(b)]{AW2013}. 

For \ref{Lemma:Existence_Double_Filtrations:Filtration}, fix a finite generating set \(m_1,\dots,m_r\) of \(M\). Let \(F_0M\) be the \(F_0A\)-submodule generated by the elements \(m_1,\dots,m_r\). This is an \(R\)-lattice by \ref{Lemma:Existence_Double_Filtrations:Lattice}. Moreover, the associated slice \(\gr_0 M=F_0M/\pi F_0M\) is a finitely generated \(\gr_0 A\)-module. Note that so far we have not made use of the chose filtration on \(\gr_0 A\). Now set \(F_i\gr_0 M=\sum_{j=1}^r\bigl(F_i\gr_0 A\bigr)\overline m_j\). This defines a filtration \(F_\bullet\gr_0 M\) for which the associated graded \(\gr\bigl(\gr_0 M\bigr)=\Gr(M)\) is a finitely generated \(\gr\bigl(\gr_0 A\bigr)=\Gr(A)\)-module (cf.\ \cite[Remark~I.5.2]{HvO1996}). Since \(\gr_0 A\) is complete with respect to the given filtration and \(\Gr(A)\) is noetherian by assumption, the filtration \(F_\bullet\gr_0 M\) is also separated by \cite[Proposition~II.2.2.1 and Theorem~I.4.14]{HvO1996}. 
\end{Proof}

\begin{Remark}\label{Remark:Good_Double_Filtrations_Independence}
Let \(A\) be complete doubly filtered \(K\)-algebra and \(M\) a finitely generated \(A\)-module as in \autoref{Lemma:Existence_Double_Filtrations}. Suppose that the associated slice \(\gr_0 A\) of \(A\) admits another filtration \(F_\bullet'\gr_0 A\) for which \(A\) is a complete doubly filtered \(K\)-algebra as in \autoref{Lemma:Existence_Double_Filtrations}. Then the good double filtrations of \(M\) constructed in \autorefpart{Lemma:Existence_Double_Filtrations}{Lemma:Existence_Double_Filtrations:Filtration} relative to either double filtration on \(A\) have the same underlying \(R\)-lattice \(F_0M\) and slice \(\gr_0 M\).
\end{Remark}

We obtain the following as an analogue and consequence of \autoref{Proposition:Stabilisation_Deformable_Algebra}.

\begin{Proposition}\label{Proposition:Stabilisation_Doubly_Filtered_Algebra}
Let \(U\) be a deformable \(R\)-algebra satisfying \eqref{Assumption} and let \(M\) be a \(\pi\)-torsionfree finitely generated \(U\)-module. Then there exists an integer \(n_0\ge0\) so that for all \(n\ge n_0\), the finitely generated \(\widehat{U_{n,K}}\)-module \(\widehat{U_{n,K}}\otimes_U M\) admits a good double filtration so that the natural morphism
\[
\Gr\bigl(\widehat{U_{n_0,K}}\otimes_U M\bigr)
\longto
\Gr\bigl(\widehat{U_{n,K}}\otimes_U M\bigr)
\]
is a graded \(K\)-linear isomorphism. 

Moreover, this isomorphism identifies the \(\Gr\bigl(\widehat{U_{n_0,K}}\bigr)\)-module structure on \(\Gr\bigl(\widehat{U_{n_0,K}}\otimes_U M\bigr)\) with the \(\Gr\bigl(\widehat{U_{n,K}}\bigr)\)-module structure on \(\Gr\bigl(\widehat{U_{n,K}}\otimes_U M\bigr)\) along the isomorphism \(\Gr\bigl(\widehat{U_{n_0,K}}\otimes_U M\bigr)\cong\Gr\bigl(\widehat{U_{n,K}}\otimes_U M\bigr)\) of \autoref{Remark:Iterated_Double_Filtrations}.
\end{Proposition}

\begin{Proof}
We reduce to \autoref{Proposition:Stabilisation_Deformable_Algebra}. Let \(n_0\ge0\) be an integer as in \autoref{Lemma:Torsion_Deformed_Module} and fix a good filtration \(F_\bullet M\) on \(M\) to define \(M_n\).

Recall that \(\widehat{M_n}\) denotes the \(\pi\)-adic completion of the \(R\)-module \(M_n\) and that \(M_K\) denotes the extension of scalars from \(R\) to \(K\) of the \(R\)-module \(M\). As \(M\) is a \(\pi\)-torsionfree finitely generated \(U\)-module, it embeds into the finitely generated \(U_K\)-module \(M_K\). The proof of \cite[Theorem~2.13]{Bode2026} now shows that the natural morphism
\[
\widehat{U_{n,K}}\otimes_{U_K} M_K
\longto\widehat{M_{n,K}}
\]
is an isomorphism for all \(n\ge n_0\). As \(\widehat{U_{n,K}}\otimes_{U_K} M_K\cong \widehat{U_{n,K}}\otimes_U M\), this yields a natural isomorphism
\[
\widehat{U_{n,K}}\otimes_U M
\longto\widehat{M_{n,K}}
\]
for all \(n\ge n_0\). 

Note that the image of \(\widehat{M_n}\) defines an \(R\)-lattice in \(\widehat{M_{n,K}}\). Indeed, by construction, we know that \(\widehat{M_n}\) is a finitely generated \(\widehat{U_n}\)-module. Since \(M_n\) is \(\pi\)-torsionfree and the \(\pi\)-multiplication map on \(M_n\) is strict for the \(\pi\)-adic filtration, we conclude that \(\widehat{M_n}\) is \(\pi\)-torsionfree from \cite[Theorem~I.3.4.13]{HvO1996}. As \(F_0\widehat{U_{n,K}}=\widehat{U_n}\), we then see that \(\widehat{M_n}\) is an \(R\)-lattice in \(\widehat{M_{n,K}}\) by \autorefpart{Lemma:Existence_Double_Filtrations}{Lemma:Existence_Double_Filtrations:Lattice}. The corresponding slice is \(\gr_0\widehat{M_{n,K}}=\widehat{M_n}/\pi\widehat{M_n}\cong M_n/\pi M_n\). The natural good filtration of \(M_n\) as a finitely generated \(U_n\)-module is separated, as \(U_n\) is positively filtered. Thus, we have defined a good double filtration on \(\widehat{U_{n,K}}\otimes_U M\cong\widehat{M_{n,K}}\) so that \(\Gr\bigl(\widehat{U_{n,K}}\otimes_U M\bigr)\cong\gr\bigl(M_n/\pi M_n\bigr)\) for all \(n\ge n_0\).

As \(\gr M_n\) is \(\pi\)-torsionfree, we find that \(\gr\bigl(M_n/\pi M_n\bigr)\cong\gr M_n/\pi\gr M_n\) as in \autoref{Proposition:Graded_Isomorphism_Doubly_Filtered_Algebra}. But by \autoref{Proposition:Stabilisation_Deformable_Algebra}, the natural morphism
\[
\mu_{n_0,n}\colon\gr M_{n_0}\longto\gr M_n
\]
is a graded \(R\)-linear isomorphism for all \(n\ge n_0\), identifying the \(\gr U_{n_0}\)-module structure \(\gr M_{n_0}\) with the \(\gr U_n\)-module structure on \(\gr M_n\) along the isomorphism \(\xi_{n_0,n}\colon\gr U_{n_0}\to\gr U_n\) for all \(n\ge n_0\). Thus, using \autoref{Proposition:Graded_Isomorphism_Doubly_Filtered_Algebra} and \autoref{Remark:Iterated_Double_Filtrations}, we conclude that for all \(n\ge n_0\) the induced natural morphism
\[
\Gr\bigl(\widehat{U_{n_0,K}}\otimes_U M\bigr)
\longto
\Gr\bigl(\widehat{U_{n,K}}\otimes_U M\bigr)
\]
is a graded \(K\)-linear isomorphism, identifying the \(\Gr\bigl(\widehat{U_{n_0,K}}\bigr)\)-module structure on \(\Gr\bigl(\widehat{U_{n_0,K}}\otimes_U M\bigr)\) with the \(\Gr\bigl(\widehat{U_{n,K}}\bigr)\)-module structure on \(\Gr\bigl(\widehat{U_{n,K}}\otimes_U M\bigr)\) along the isomorphism \(\Gr\bigl(\widehat{U_{n_0,K}}\otimes_U M\bigr)\cong \Gr\bigl(\widehat{U_{n,K}}\otimes_U M\bigr)\).
\end{Proof}

\subsection{Application to dimension theory of doubly filtered \texorpdfstring{\(K\)}{K}-algebras} Let \(A\) be a doubly filtered \(K\)-algebra such that \(\Gr(A)\) is commutative and noetherian and let \(M\) be a finitely generated \(A\)-module. Choose a good double filtration on \(M\). The \bemph{characteristic variety of \(\bm M\)} is the subset
\[
\bm{\Ch(M)}
\coloneqq\Supp\bigl(\Gr(M)\bigr)
\subseteq\Spec\bigl(\Gr(A)\bigr)
\]

\begin{Remark}\label{Remark:Support_Annihilator}
For a good double filtration of \(M\), the associated graded \(\Gr(M)\) is a finitely generated \(\Gr(A)\)-module by definition. Thus, \(\Supp\bigl(\Gr(M)\bigr)=V\bigl(\sqrt{\Ann \Gr(M)}\bigr)\) and the characteristic variety is a closed subset. 
\end{Remark}

The characteristic variety is an intrinsic invariant of \(M\) and behaves as expected in short exact sequences.

\begin{Proposition}[{{\cite[Proposition~3.3]{AW2013}}}]\label{Proposition:Basic_Properties_Characteristic_Variety}
Let \(A\) be a doubly filtered \(K\)-algebra such that \(\Gr(A)\) is commutative and noetherian. For any finitely generated \(A\)-module \(M\), the characteristic variety \(\Ch(M)\) does not depend on the choice of a good double filtration. Moreover, if
\[
\begin{tikzcd}
0 \ar[r] & M' \ar[r] & M \ar[r] & M'' \ar[r] & 0    
\end{tikzcd}
\]
is a short exact sequence of finitely generated \(A\)-modules, then there are good double filtrations on \(M'\), \(M\) and \(M''\) such that there is a short exact sequence
\[
\begin{tikzcd}
0 \ar[r] & \Gr(M') \ar[r] & \Gr(M) \ar[r] & \Gr(M'') \ar[r] & 0    
\end{tikzcd}
\]
of finitely generated \(\Gr(A)\)-modules. In particular, \(\Ch(M)=\Ch(M')\cup\Ch(M'')\).
\end{Proposition}

Let \(A\) be a doubly filtered \(K\)-algebra such that \(\Gr(A)\) is commutative and noetherian and let \(M\) be a finitely generated \(A\)-module. Choose a good double filtration on \(M\). The \bemph{dimension of \(\bm M\)} is
\[
\bm{\dim(M)}
\coloneqq\dim\Ch(M)\,.
\]
By \autoref{Proposition:Basic_Properties_Characteristic_Variety}, this is well-defined. The following is now immediate from \autoref{Proposition:Basic_Properties_Characteristic_Variety} and \autoref{Proposition:Stabilisation_Doubly_Filtered_Algebra}.

\begin{Corollary}
Let \(U\) be a deformable \(R\)-algebra satisfying \eqref{Assumption} for which \(\gr U\) is commutative and let \(M\) be a \(\pi\)-torsionfree finitely generated \(U\)-module. Then there exists an integer \(n_0\ge0\) so that
\[
\dim\bigl(\widehat{U_{n,K}}\otimes_U M\bigr)
=\dim\bigl(\widehat{U_{n_0,K}}\otimes_U M\bigr)
\]
for all \(n\ge n_0\).
\end{Corollary}

\section{On Completed Weyl Algebras}\label{Section:Completed_Weyl_Algebras}

We show that completed Weyl algebras admit two structures as doubly filtered \(K\)-algebras, using either the order filtration or the Bernstein filtration. The latter gives rise to a notion of Hilbert polynomial for finitely generated modules over completed Weyl algebras.

\subsection{Generalities on completed Weyl algebras} We start by collecting some generalities on the Weyl algebras over general commutative base rings, our cases of interests being the rings \(R\) and \(k\).

\begin{Definition}
Let \(S\) a commutative ring. The \bemph{\(\bm d\)-th Weyl algebra over \(\bm S\)}, denoted \(\bm{A_d(S)}\), is the noncommutative \(S\)-algebra with \(2d\) generators \(x_1,\dots,x_d,y_1,\dots,y_d\) and relations
\[
[x_i,x_j]=x_ix_j-x_jx_i=0
,\qquad
[y_i,y_j]=y_iy_j-y_jy_i=0
,\qquad
[y_i,x_j]=y_ix_j-x_jy_i=\delta_{ij}\,.
\]
\end{Definition}

For \(S\) a commutative ring, the \(d\)-th Weyl algebra \(A_d(S)\) can be understood as an iterated skew polynomial extension of \(S[x_1,\dots,x_d\)], as described in \cite[Chapter~2,~p.30]{GW1989} (note the interchanged role of \(x_i\) and \(y_i\)). See also \cite[Paragraph~1.3.8]{MCR2001}, where the sign convention for skew polynomial extensions is switched. In particular, using multi-index notation, every element of \(A_d(S)\) admits a standard form \(\sum_{\alpha,\beta} s_{\alpha,\beta}x^\alpha y^\beta\) with \(s_{\alpha,\beta}\in S\). From this we can define two (positive) \(S\)-algebra filtrations of \(A_d(S)\) (cf.\ \cite[Example~I.2.3.I]{HvO1996}). The \bemph{Bernstein filtration \(\bm{\cB_{\bullet}}\) on \(\bm{A_d(S)}\)} is the positive \(S\)-algebra filtration with \(\cB_i A_d(S)\) the \(S\)-submodule generated by \(x^\alpha y^\beta\) with \(|\alpha|+|\beta|\le i\). The \bemph{order filtration \(\bm{\cF_{\bullet}}\) on \(\bm{A_d(S)}\)} is the positive \(S\)-algebra filtration with \(\cF_i A_d(S)\) the \(S[x_1,\dots,x_d]\)-submodule generated by \(y^\beta\) with \(|\beta|\le i\).

\begin{Proposition}\label{Proposition:Weyl_Algebra_Associated_Graded}
Let \(S\) be a commutative ring. The associated graded rings of \(A_d(S)\) for the Bernstein filtration \(\cB_\bullet A_d(S)\) and the order filtration \(\cF_\bullet A_d(S)\) admit the following descriptions.
\begin{enumroman}
\item\label{Proposition:Weyl_Algebra_Associated_Graded:Bernstein} The graded \(S\)-algebra \(\gr^{\cB_\bullet}\bigl(A_d(S)\bigr)\) is isomorphic to the polynomial ring \(S[X_1,\dots,X_d,Y_1,\dots,Y_d]\) in \(2d\) variables over \(S\) equipped with the standard grading.
\item\label{Proposition:Weyl_Algebra_Associated_Graded:Order} The graded \(S\)-algebra \(\gr^{\cF_\bullet}\bigl(A_d(S)\bigr)\) is isomorphic to the polynomial ring \(S[X_1,\dots,X_d][Y_1,\dots,Y_d]\) in \(d\) variables over \(S[X_1,\dots,X_d]\) equipped with the standard grading.
\end{enumroman}
In particular, the graded \(S\)-algebras \(\gr^{\cB_\bullet}\bigl(A_d(S)\bigr)\) and \(\gr^{\cF_\bullet}\bigl(A_d(S)\bigr)\) are isomorphic as abstract rings to a polynomial ring in \(2d\) variables over \(S\).
\end{Proposition}

\begin{Proof}
Using the standard form introduced above, the proof of \ref{Proposition:Weyl_Algebra_Associated_Graded:Bernstein} can be carried out as in \cite[Theorem~7.3.1]{Coutinho1995}. The proof of \ref{Proposition:Weyl_Algebra_Associated_Graded:Order} is similar (cf.\ \cite[Exercise~7.6.5]{Coutinho1995}, \cite[Section~1.2]{HTT2008}).
\end{Proof}

\begin{Corollary}
The \(d\)-th Weyl algebra \(A_d(R)\) over \(R\) is a deformable \(R\)-algebra satisfying \eqref{Assumption} with respect to either the Bernstein filtration \(\cB_\bullet A_d(R)\) or the order filtration \(\cF_\bullet A_d(R)\).
\end{Corollary}

\begin{Definition}
Let \(A_d(R)\) be the \(d\)-th Weyl algebra over \(R\) equipped with the order filtration \(\cF_\bullet A_d(R)\). Denote the \(n\)-th deformation by \(\bm{A_{d,n}(R)}\coloneqq A_d(R)_n\). We call \(\bm{\widehat{A_{d,n}(K)}}\coloneqq\widehat{A_{d,n}(R)_K}\) the \bemph{\(\bm n\)-th deformed completed \(\bm d\)-th Weyl algebra} and \(\bm{\widehat{A_d(K)}}\coloneqq\widehat{A_{d,0}(K)}\) the \bemph{completed \(\bm d\)-th Weyl algebra}.
\end{Definition}

We will collectively refer to the \(n\)-th deformed completed (\(d\)-th) Weyl algebras as completed (\(d\)-th) Weyl algebras. The completed Weyl algebras are doubly filtered \(K\)-algebras by \autoref{Proposition:Graded_Isomorphism_Doubly_Filtered_Algebra}. Their canonical slices admit the following description.

\begin{Lemma}\label{Lemma:Completed_Weyl_Algebras_Slices}
\begin{enumroman}
\item The canonical slice of the doubly filtered \(K\)-algebra \(\widehat{A_d(K)}\) is isomorphic to the \(d\)-th Weyl algebra \(A_d(k)\) over the residue field \(k\).
\item If \(n>0\), then the canonical slice of the doubly filtered \(K\)-algebra \(\widehat{A_{d,n}(K)}\) is isomorphic to the polynomial ring \(k[X_1,\dots,X_d,Y_1,\dots,Y_d]\) over the residue field \(k\).
\end{enumroman}
\end{Lemma}

\begin{Proof}
The proof of \autoref{Proposition:Graded_Isomorphism_Doubly_Filtered_Algebra} (cf.\ \cite[Lemma~3.7]{AW2013}) shows that the doubly filtered \(K\)-algebra structure on \(\widehat{A_{d,n}(K)}=\widehat{A_d(R)_{n,K}}\) is obtained by considering the \(R\)-lattice \(\widehat{A_d(R)_n}\) and endowing the associated slice \(\gr_0\widehat{A_{d,n}(K)}=\widehat{A_d(R)_n}/\pi\widehat{A_d(R)_n}\) with the induced filtration. Note that there is a canonical isomorphism \(\widehat{A_d(R)_n}/\pi\widehat{A_d(R)_n}\cong A_d(R)_n/\pi A_d(R)_n\) since we consider \(\pi\)-adic completions.

For the order filtration \(\cF_\bullet A_d(R)\), the \(n\)-th deformation is the \(R\)-subalgebra generated by \(2d\) generators \(x_1,\dots,x_d,\pi^n y_1,\dots,\pi^ny_d\). To see this, recall that the induced filtration \(\cF_\bullet A_d(R)_n\) is explicitly given by \(\cF_i A_d(R)_n=\sum_{j=0}^i\pi^{jn}\cF_j A_d(R)\). Thus, \(A_{d,n}(R)\) is isomorphic to the \(R\)-algebra with \(2d\) generators \(\tilde x_1,\dots,\tilde x_d,\tilde y_1,\dots,\tilde y_d\) subject to the relations
\[
[\tilde x_i,\tilde x_j]=0
,\qquad
[\tilde y_i,\tilde y_j]=0
,\qquad
[\tilde y_i,\tilde x_j]=\pi^n\delta_{ij}\,.
\]
In particular, reducing modulo \(\pi\) we recover the presentation of \(A_d(k)\) for \(n=0\), while we recover the presentation of \(k[X_1,\dots,X_d,Y_1,\dots,Y_d]\) for \(n>0\), where \(\tilde x_i\) and \(\tilde y_i\) correspond to \(X_i\) and \(Y_i\), respectively (cf.\ \cite[Section~5.2]{Bode2025Almost}). 
\end{Proof}

\begin{Corollary}\label{Corollary:Associated_Graded_Order}
There is an isomorphism
\[
\Gr^{\cF_\bullet}\bigl(\widehat{A_{d,n}(K)}\bigr)
\cong k[X_1,\dots,X_d][Y_1,\dots,Y_d]
\]
of graded \(k\)-algebras for all \(n\ge0\).
\end{Corollary}

\begin{Proof}
By \autoref{Proposition:Graded_Isomorphism_Doubly_Filtered_Algebra} and \autoref{Lemma:Completed_Weyl_Algebras_Slices} we know that
\[
\Gr^{\cF_\bullet}\bigl(\widehat{A_{d,n}(K)}\bigr)
\cong\gr^{\cF_\bullet}\bigl(\gr_0\widehat{A_{d,n}(K)}\bigr)
\cong
\begin{cases}
\gr^{\cF_\bullet}\bigl(A_d(k)\bigr) & \text{if}\ n=0
\\
\gr^{\cF_\bullet}\bigl(k[X_1,\dots,X_d][Y_1,\dots,Y_d]\bigr) & \text{if}\ n>0
\end{cases}
\]
where the slices carry the filtration induced from \(A_d(R)_n\). This is the order filtration on \(A_d(k)\) for \(n=0\) and the degree filtration on \(k[X_1,\dots,X_d][Y_1,\dots,Y_d]\) for \(n>0\). In particular, it follows that
\[
\Gr^{\cF_\bullet}\bigl(\widehat{A_{d,n}(K)}\bigr)
\cong k[X_1,\dots,X_d][Y_1,\dots,Y_d]
\]
for all \(n\ge0\). For \(n=0\), we use \autorefpart{Proposition:Weyl_Algebra_Associated_Graded}{Proposition:Weyl_Algebra_Associated_Graded:Order} while the case \(n>0\) is clear.
\end{Proof}

\begin{Corollary}\label{Corollary:Noetherian_Auslander_Regular}
For each \(n\ge0\), the complete doubly filtered \(K\)-algebra \(\widehat{A_{d,n}(K)}\) is noetherian and Auslander regular of global dimension at most \(2d\).
\end{Corollary}

\begin{Proof}
Noetherianity and Auslander regularity follow from \cite[Theorem~3.3]{AW2013} (loc.\ cit.\ Auslander regularity includes noetherianity). The assertion on the global dimension follows from the proof of \cite[Theorem~3.3]{AW2013} combined with \cite[Theorem~II.3.1.4(2)]{HvO1996} and \(\dim k[X_1,\dots,X_d][Y_1,\dots,Y_d]=2d\).
\end{Proof}

We shall also need a different double filtration on the completed Weyl algebras. Given our explicit description of the slices, this is immediate.

\begin{Corollary}\label{Corollary:Associated_Graded_Bernstein}
For each \(n\ge0\), the canonical slice \(\gr_0\widehat{A_{d,n}(K)}\) of the \(R\)-lattice \(\widehat{A_{d,n}(R)}\) admits a filtration \(\cB_\bullet\gr_0\widehat{A_{d,n}(K)}\) such that \(\widehat{A_{d,n}(K)}\) is a complete doubly filtered \(K\)-algebra with associated graded \(k\)-algebra
\[
\Gr^{\cB_\bullet}\bigl(\widehat{A_{d,n}(K)}\bigr)
\cong k[X_1,\dots,X_d,Y_1,\dots,Y_d]\,.
\]
\end{Corollary}

\begin{Proof}
From \autoref{Lemma:Completed_Weyl_Algebras_Slices} we know that the slices of the \(R\)-lattices \(\widehat{A_{d,n}(R)}\) are
\[
\gr_0\widehat{A_{d,n}(K)}
\cong
\begin{cases}
A_d(k) & \text{if}\ n=0
\\
k[X_1,\dots,X_d,Y_1,\dots,Y_d] & \text{if}\ n>0\,.
\end{cases}
\]
For \(n=0\) choose the Bernstein filtration and for \(n>0\) choose the standard degree filtration. As these filtrations are positive, they are complete. Thus, they define a complete doubly filtered \(K\)-algebra structure on \(\widehat{A_{d,n}(K)}\). Moreover, we find
\[
\Gr^{\cB_\bullet}\bigl(\widehat{A_{d,n}(K)}\bigr)
\cong k[X_1,\dots,X_d,Y_1,\dots,Y_d]
\]
for all \(n\ge0\). For \(n=0\), we use \autorefpart{Proposition:Weyl_Algebra_Associated_Graded}{Proposition:Weyl_Algebra_Associated_Graded:Bernstein} while the case \(n>0\) is clear.
\end{Proof}

\subsection{Dimension theory and Hilbert polynomials for completed Weyl algebras}\label{Subsection:Dimensions_Completed_Weyl_Algebras} In particular, \autoref{Corollary:Associated_Graded_Order} and \autoref{Corollary:Associated_Graded_Bernstein} show that the complete double filtrations we consider have as associated graded a commutative noetherian ring. Thus, we can define characteristic varieties for finitely generated modules. Crucially, the dimensions we obtain with respect to either complete double filtration agree!

In fact, both dimensions agree with a purely homological notion of dimension. Recall that for a ring \(A\) and an \(A\)-module \(M\), the \bemph{homological grade of \(\bm M\)} is defined to be \(j_A(M)\coloneqq\min\bigl\{i\,|\,\Ext_A^i(M,A)\ne0\bigr\}\) or \(\infty\) if no such \(i\) exists.

\begin{Proposition}\label{Proposition:Dimensions_Equal}
Let \(n\ge0\) and let \(M_n\) be a finitely generated \(\widehat{A_{d,n}(K)}\)-module. Then
\[
\dim\bigl(\Ch^{\cF_\bullet}(M_n)\bigr)
=2d-j_{\widehat{A_{d,n}(K)}}(M_n)
=\dim\bigl(\Ch^{\cB_\bullet}(M_n)\bigr)\,.
\]
\end{Proposition}

\begin{Proof}
Recall that there are isomorphisms of \(k\)-algebras
\[
\Gr^{\cF_\bullet}\bigl(\widehat{A_{d,n}(K)}\bigr)
\cong k[X_1,\dots,X_d,Y_1,\dots,Y_d]
\cong\Gr^{\cB_\bullet}\bigl(\widehat{A_{d,n}(K)}\bigr)
\]
by \autoref{Corollary:Associated_Graded_Order} and \autoref{Corollary:Associated_Graded_Bernstein}. This identification links the two filtrations. 

By construction, both complete double filtrations on \(\widehat{A_{d,n}(K)}\) share the lattice \(\widehat{A_{d,n}(R)}\) and the slice \(\gr_0\widehat{A_{d,n}(K)}\). Moreover, both their associated graded \(k\)-algebras are noetherian. Thus, as discussed in \autoref{Remark:Good_Double_Filtrations_Independence}, the good double filtrations of the finitely generated module \(M_n\) constructed in \autorefpart{Lemma:Existence_Double_Filtrations}{Lemma:Existence_Double_Filtrations:Filtration} share the same finitely generated \(\gr_0\widehat{A_{d,n}(K)}\)-module \(\gr_0 M_n\) as a slice. Moreover, by \autoref{Proposition:Basic_Properties_Characteristic_Variety} we may use any good double filtration to compute the characteristic variety.

Now, observe that \(k[X_1,\dots,X_d,Y_1,\dots,Y_d]\) is a polynomial ring over a field \(k\), hence a regular commutative noetherian ring of pure dimension \(2d\). Applying \cite[Theorem~D.4.3]{HTT2008} for both filtrations, we see that
\[
\dim\biggl(\Supp\bigl(\gr^{\cF_\bullet}(\gr_0 M_n)\bigr)\biggr)
=2d-j_{\gr_0\widehat{A_{d,n}(K)}}(\gr_0 M_n)
=\dim\biggl(\Supp\bigl(\gr^{\cB_\bullet}(\gr_0 M_n)\bigr)\biggr)\,.
\]
In order to conclude, we check that \(j_{\gr_0\widehat{A_{d,n}(K)}}(\gr_0 M_n)=j_{\widehat{A_{d,n}(K)}}(M_n)\). Denote the associated graded of the \(\pi\)-adic filtration on \(\widehat{A_{d,n}(K)}\) determined by the \(R\)-lattice \(\widehat{A_{d,n}(R)}\) by \(\gr^\pi\widehat{A_{d,n}(K)}\) and denote by \(\gr^\pi M_n\) the associated graded of the filtration on \(M_n\) determined by our choice of generators. By \autoref{Corollary:Noetherian_Auslander_Regular}, we know that \(\gr^\pi\widehat{A_{d,n}(K)}\) is noetherian and Auslander regular. Moreover, \(\gr^\pi\widehat{A_{d,n}(K)}\cong\bigl(\gr_0\gr^\pi\widehat{A_{d,n}(K)}\bigr)[t^{\pm1}]\) by \cite[Lemma~3.1]{AW2013}. Since \(\gr^\pi\widehat{A_{d,n}(K)}\) is also strongly \(\bZ\)-graded in the sense of \cite[Section~I.4.1]{HvO1996}, we conclude that \(j_{\gr^\pi\widehat{A_{d,n}(K)}}(\gr^\pi M_n)=j_{\gr_0\widehat{A_{d,n}(K)}}(\gr_0 M_n)\) from \cite[Lemma~III.2.5.3]{HvO1996}. Thus, \(j_{\widehat{A_{d,n}(K)}}(M_n)=j_{\gr^\pi\widehat{A_{d,n}(K)}}(\gr^\pi M_n)\) by \cite[Theorem~III.2.5.2]{HvO1996}.
\end{Proof}

For the order filtration, \cite[Corollary~7.4]{AW2013} shows the analogue of Bernstein's inequality, meaning that for \(n\ge0\) and any non-zero finitely generated \(\widehat{A_{d,n}(K)}\)-module \(M_n\), we have \(\dim\bigl(\Ch^{\cF_\bullet}(M_n)\bigr)\ge d\). By \autoref{Proposition:Dimensions_Equal} we obtain the following.

\begin{Corollary}\label{Corollary:Bernstein_Inequality_Bernstein}
Let \(n\ge0\) and let \(M_n\) be a non-zero finitely generated \(\widehat{A_{d,n}(K)}\)-module. Then 
\[
\dim\bigl(\Ch^{\cB_\bullet}(M_n)\bigr)\ge d\,.
\]
\end{Corollary}

The main reason why we have introduced the double filtration on \(\widehat{A_{d,n}(K)}\) based on the Bernstein filtration is to make use of the formalism of Hilbert polynomials (see \autoref{Appendix:Hilbert_Polynomials}). Since the \(k\)-algebra \(\gr_0\widehat{A_{d,n}(K)}\) is finitely generated by \(2d\) generators and the filtration \(\cB_\bullet\gr_0\widehat{A_{d,n}(K)}\) is standard finite-dimensional  with \(\Gr^{\cB_\bullet}\bigl(\widehat{A_{d,n}(K)}\bigr)\) a polynomial ring in \(2d\) variables over \(k\) by \autoref{Corollary:Associated_Graded_Bernstein}, the \(k\)-algebra \(\gr_0\widehat{A_{d,n}(K)}\) is almost commutative in the sense of \autoref{Definition:Somewhat_Commutative}. Hence, we can use \autoref{Definition:Hilbert_Polynomial} and \autoref{Definition:Dimension_Multiplicity}.

\begin{Definition}\label{Definition:Hilbert_Polynomial_Completed_Weyl_Algebras}
Let \(n\ge0\) and let \(M_n\) be a finitely generated \(\widehat{A_{d,n}(K)}\)-module. Let \((F_0 M_n,\cB_\bullet\gr_0 M_n)\) be a good double filtration relative to the complete double filtration \((\widehat{A_{d,n}(R)},\cB_\bullet\widehat{A_{d,n}(R)})\). A \bemph{Hilbert polynomial of \(\bm {M_n}\)} is a Hilbert polynomial of the slice \(\gr_0 M_n\). The \bemph{dimension of \(\bm{M_n}\)}, denoted by \(\bm{\rd_n(M_n)}\), and the \bemph{multiplicity of \(\bm{M_n}\)}, denoted by \(\bm{\rm_n(M_n)}\), are the dimension and the multiplicity of \(\gr_0 M_n\).
\end{Definition}

We proceed as in the proof of \autoref{Proposition:Basic_Properties_Characteristic_Variety} to see that the dimension and multiplicity are independent of the choice of a good double filtration.

\begin{Proposition}
Let \(n\ge0\) and let \(M_n\) be a finitely generated \(\widehat{A_{d,n}(K)}\)-module. Then \(\rd_n(M_n)\) and \(\rm_n(M_n)\) are independent of the choice of a good double filtration.
\end{Proposition}

\begin{Proof}
Let \((F_0 M_n,\cB_\bullet\gr_0 M_n)\) and \((F_0' M_n,\cB_\bullet\gr_0' M_n)\) be good double filtrations of \(M_n\). We now proceed as in the proof of \cite[Proposition~1.1.2]{Ginsburg1986}. First, suppose that the lattices \(F_0M_n\) and \(F_0'M_n\) satisfy \(\pi\cdot F_0'M_n\subseteq \pi\cdot F_0M_n\subseteq F_0'M_n\subseteq F_0 M_n\). We will call such lattices neighbouring. Then there are short exact sequences
\[
\begin{tikzcd}[row sep=small]
0 \ar[r] &  F_0'M_n/\pi\cdot F_0M_n \ar[r] & \gr_0 M_n \ar[r] & F_0M_n/F_0'M_n \ar[r] & 0
\\
0 \ar[r] &  \pi\cdot F_0M_n/\pi\cdot F_0'M_n \ar[r] & \gr_0' M_n \ar[r] & F_0'M_n/\pi\cdot F_0 M_n \ar[r] & 0
\end{tikzcd}
\]
of finitely generated \(\gr_0\widehat{A_{d,n}(K)}\)-modules. Now \autorefpart{Proposition:Dimension_Multiplicity_SES}{Proposition:Dimension_Multiplicity_SES:Dimension} shows that
\begin{align*}
\rd(\gr_0 M_n)
&=\max\big(\rd(F_0'M_n/\pi\cdot F_0M_n),\rd(F_0M_n/F_0'M_n)\big)
\\
\rd(\gr_0' M_n)
&=\max\big(\rd(\pi\cdot F_0M_n/\pi\cdot F_0'M_n),\rd(F_0'M_n/\pi\cdot F_0 M_n)\big)
\end{align*}
Using that \(\pi\cdot F_0M_n/\pi\cdot F_0'M_n\cong F_0M_n/F_0'M_n\), this shows \(\rd(\gr_0 M_n)=\rd(\gr_0' M_n)\) for neighbouring lattices. For each of the cases in \autorefpart{Proposition:Dimension_Multiplicity_SES}{Proposition:Dimension_Multiplicity_SES:Cases}, we also conclude that \(\rm(\gr_0 M_n)=\rm(\gr_0' M_n)\). This proves the claim for neighbouring lattices.

By \cite[Lemma~1.1.1]{Ginsburg1986}, there are integers \(k,l\ge0\) such that \(\pi^k\cdot F_0 M_n\subseteq F_0' M_n\subseteq\pi^{-l}\cdot F_0 M_n\). Now consider the lattices \((F_0 M)_j=F_0 M+\pi^j F_0'M\) for \(j\in\bZ\). We see that \((F_0 M)_j=F_0 M\) for \(j\gg0\) and \((F_0 M)_j=\pi^j F_0' M\) for \(j\ll0\). Moreover, the lattices \((F_0 M)_j\) and \((F_0 M)_{j+1}\) are neighbouring for all \(j\in\bZ\) and \(\pi^j F_0'M/\pi^{j+1} F_0'M\cong\gr_0' M\). In particular, \(\rd(\gr_0 M_n)=\rd(\gr_0' M_n)\) and \(\rm(\gr_0 M_n)=\rm(\gr_0' M_n)\), which implies the claim.
\end{Proof}

The notion of dimension just introduced agrees with our earlier considerations.

\begin{Proposition}
Let \(n\ge0\) and let \(M_n\) be a finitely generated \(\widehat{A_{d,n}(K)}\)-module. Then
\[
\rd_n(M_n)
=\dim\bigl(\Ch^{\cB_\bullet}(M_n)\bigr)
=\dim\bigl(\Ch^{\cF_\bullet}(M_n)\bigr)\,.
\]
\end{Proposition}

\begin{Proof}
Using \autoref{Proposition:Hilbert_Polynomial_Degree_Dimension}, we have
\[
\rd_n(M_n)
=\rd\big(\gr_0 M_n\big)
=\dim V\bigl(\Ann\gr^{\cB_\bullet}(\gr_0 M_n)\bigr)
=\dim V\bigl(\Ann\Gr^{\cB_\bullet}(M_n)\bigr)\,.
\]
As 
\[
V\bigl(\Ann\Gr^{\cB_\bullet}(M_n)\bigr)
=V\biggl(\sqrt{\Ann\Gr^{\cB_\bullet}(M_n)}\biggr)
=\Supp\bigl(\Gr^{\cB_\bullet}(M_n)\bigr)
=\Ch^{\cB_\bullet}(M_n)
\]
by \autoref{Remark:Support_Annihilator}, the first equality follows. The second equality is \autoref{Proposition:Dimensions_Equal}.
\end{Proof}

To close this section, we introduce holonomicity for completed Weyl algebras and adapt results from the classical setting, using the dimension and multiplicity we just defined in \autoref{Definition:Hilbert_Polynomial_Completed_Weyl_Algebras}

\begin{Definition}\label{Definition:Holonomic_Completed_Weyl_Algebra}
Let \(n\ge0\) and let \(M_n\) be a finitely generated \(\widehat{A_{d,n}(K)}\)-module. We say that \(M_n\) is \bemph{holonomic} if \(M_n=0\) or \(\rd_n(M_n)=d\).
\end{Definition}

We have the following properties of the dimension and the multiplicities associated to finitely generated modules. For the classical setting, cf. \cite[Section~9]{Coutinho1995}. Using \autoref{Proposition:Basic_Properties_Characteristic_Variety}, their proof goes through as in the \cite[Section~8.6]{MCR2001} (see also \autoref{Proposition:Dimension_Multiplicity_SES}).

\begin{Proposition}\label{Proposition:Properties_Dimension_Multiplicity}
Let \(n\ge0\) and let
\[
\begin{tikzcd}
0 \ar[r] & M_n' \ar[r] & M_n \ar[r] & M_n'' \ar[r] & 0
\end{tikzcd}
\]
be short exact sequence of finitely generated \(\widehat{A_{d,n}(K)}\)-modules.
\begin{enumroman}
\item\label{Proposition:Properties_Dimension_Multiplicity:Multiplicity} If \(M_n\) is non-zero, then \(\rm_n(M_n)>0\).
\item\label{Proposition:Properties_Dimension_Multiplicity:Dimension_SES} \(\rd_n(M_n)=\max(\rd_n(M_n'),\rd_n(M_n''))\).
\item\label{Proposition:Properties_Dimension_Multiplicity:Multiplicity_SES} If \(\rd_n(M_n')=\rd_n(M_n)=\rd_n(M_n'')\), then \(\rm_n(M_n)=\rm_n(M_n')+\rm_n(M_n'')\).
\end{enumroman}
\end{Proposition}

Using \autoref{Proposition:Properties_Dimension_Multiplicity}, we then obtain the following central result as a direct adaptation of \cite[Theorem~10.2.2 and Scholium~10.2.3]{Coutinho1995} (see also \cite[Corollary~8.5.7]{MCR2001}). For the convenience of the reader, we sketch the proof.

\begin{Corollary}\label{Corollary:Multiplicity_Length_Bound}
Let \(n\ge0\) and let \(M_n\) be a non-zero holonomic \(\widehat{A_{d,n}(K)}\)-module. Then \(M_n\) is of finite length bounded by \(\rm_n(M_n)\).
\end{Corollary}

\begin{Proof}
Let \(M_n=N_0\gneq N_1\gneq N_2\gneq\cdots\gneq N_r\) be a strict chain of \(\widehat{A_{d,n}(K)}\)-submodules of \(M_n\). Combining \autoref{Corollary:Bernstein_Inequality_Bernstein} with \autoref{Proposition:Properties_Dimension_Multiplicity} shows that \(\rm_n(M_n)=\sum_{i=0}^{r-1}\rm_n(N_i/N_{i+1})+\rm_n(N_r)\). As we assume each inclusion to be strict, \autoref{Proposition:Properties_Dimension_Multiplicity} shows that the right-hand side is bounded from below by \(r\). In particular, any strict chain of submodules of \(M\) has length bounded by \(\rm_n(M_n)\). This shows that \(M_n\) is of finite length bounded by \(\rm_n(M_n)\).
\end{Proof}

\section{Sheaves of Differential Operators and Their Modules}\label{Section:Differential_Operators}

In this section we collect generalities on sheaves of modules over sheaves of differential operators on rigid analytic \(K\)-spaces. Differential operators are typically only considered for smooth rigid analytic \(K\)-spaces. Locally, one considers affinoid subdomains admitting coordinate systems (see \autoref{Appendix:D_Modules_Rigid_Spaces:Coordinates}). We choose to introduce both algebraic \(\cD\)-modules as well as analytic \(\Dcap\)-module using the language of Lie--Rinehart algebras, the algebraic equivalent of Lie algebroids. To keep the exposition brief, we only survey the relevant aspects of the respective theory. An exception are coherent \(\cD\)-modules, to which we devote \autoref{Appendix:D_Modules_Rigid_Spaces:Coherent} as there are some subtleties in the rigid analytic setting.

\subsection{Lie--Rinehart algebras}

Let \(S\) be a commutative ring and let \(B\) be a commutative \(S\)-algebra. Recall from \cite[Section~2]{Rinehart1963} that an \bemph{\(\bm{(S,B)}\)-Lie algebra}, is a \(B\)-module \(L\) equipped with an \(S\)-bilinear Lie bracket and a \(B\)-linear \(S\)-Lie algebra morphism
\[
\rho\colon L\longto\Der_S(B)
\]
satisfying \([x,by]=b[x,y]+\rho(x)(b)y\) for all \(b\in B\) and \(x,y\in L\). Adopting the terminology of \cite[Section~2]{AW2019}, we may also refer to \((S,B)\)-Lie algebras as \bemph{Lie--Rinehart algebra}. 

To each \((S,B)\)-Lie algebra \(L\), there is an associated \bemph{enveloping algebra \(\bm{\cU_B(L)}\)} as defined in \cite[Section~2]{Rinehart1963}. This an associative \(S\)-algebra equipped with canonical morphisms \(i_B\colon B\longto\cU_B(L)\) and \(i_L\colon L\longto U_B(L)\) of \(S\)-algebras and \(S\)-Lie algebras, respectively, which are the universal morphism satisfying \(i_L(bx)=i_B(b)i_L(x)\) and \([i_L(x),i_B(b)]=i_B(\rho(x)(b))\) for all \(b\in B\) and \(x\in L\). Note that \(U_B(L)\) has a natural degree filtration, obtained by setting 
\[
F_0\cU_B(L)=B
,\qquad
F_1\cU_B(L)=B+L
,\qquad
F_i\cU_B(L)=F_1\cU_B(L)\cdot F_{i-1}\cU_B(L)\quad \text{for}\quad i\ge2\,.
\]
We say that an \((S,B)\)-Lie algebra is \bemph{smooth} if it is coherent and projective as \(B\)-module. For a smooth \((S,B)\)-Lie algebra \(L\), there is a natural isomorphism
\[
\Sym_B L\longto \gr\cU_B(L)
\]
of \(B\)-algebras by \cite[Theorem~3.1]{Rinehart1963}.

Let \(A\) be an affinoid \(K\)-algebra. Recall that an \bemph{affine formal model for \(\bm A\)} is a topologically finitely presented \(R\)-subalgebra \(\cA\subset A\) so that \(K\otimes_R\cA=A\). Following \cite[Section~6]{AW2019}, for \(L\) a \((K,A)\)-Lie algebra and \(\cA\) an affine formal model for \(A\), we call an \((R,\cA)\)-Lie subalgebra an \bemph{\(\bm{(R,\cA)}\)-Lie lattice} if \(K\otimes_R\cL=L\) and \(\cL\) is finitely generated as an \(\cA\)-module. Observe that for any smooth \((R,\cA)\)-Lie lattice, the enveloping algebra \(\cU_\cA(\cL)\) is a deformable \(R\)-algebra. Its deformations are given by the enveloping algebras \(\cU_\cA(\pi^n\cL)\) associated to the scaled \((R,\cA)\)-Lie lattices \(\pi^n\cL\).

\subsection{Overview of algebraic \texorpdfstring{\(\cD\)}{D}-modules} We will make use of the theory of algebraic \(\cD\)-modules on rigid analytic \(K\)-space. Many considerations involving coherent algebraic \(\cD\)-modules can be carried out as in \cite{HTT2008} or \cite{Mebkhout1989}. This is explicitly done in \cite[Section~4.3]{MNM1991} and \cite[Section~2]{Raczka2024Holonomic}. While the latter uses the language of adic spaces, the results relevant to our applications translate without difficulties.

Contrary to the references given, we choose to introduce the sheaf of differential operators \(\cD\) using enveloping algebras of Lie--Rinehart algebras, as in \cite{BB2021}. If \(\cT_X\) is the tangent sheaf of a rigid analytic \(K\)-space \(X\) (see \autoref{Appendix:D_Modules_Rigid_Spaces:Coordinates}), then for any affinoid subdomain \(U\) of \(X\), the \(\cO_X(U)\)-module \(\cT_X(U)\) is an \((K,\cO_X(U))\)-Lie algebra.

\begin{Definition}[{{\cite[Definition~2.2]{BB2021}}}]\label{Definition:Algebraic_Differential_Operators_D}
Let \(X\) be a rigid analytic \(K\)-space. The \bemph{sheaf of algebraic differential operators \(\bm{\cD_X}\)} is the sheaf of \(K\)-algebras defined by
\[
\cD_X(U)
\coloneqq\cU_{\cO_X(U)}\bigl(\cT_X(U)\bigr)
\]
for \(U\) an affinoid subdomain of \(X\).
\end{Definition}

\begin{Remark}
That \(\cD_X\) is a sheaf can be seen as follows. Let \(X=\Sp A\) be an affinoid \(K\)-space and let \(L=\cT_X(X)\). Consider the \(A\)-module \(\cU_A(L)\). By \cite[Corollary~8.2.5]{BGR1984}, the presheaf determined by \(U\mapsto\cO_X(U)\otimes_A \cU_A(L)\) for affinoid subdomains \(U\) of \(X\) is a sheaf. Since for any affinoid subdomain \(U\) the restriction morphism \(A\to\cO_X(U)\) is étale, \cite[Lemma~2.4, Proposition~2.3]{AW2019} show there is a natural isomorphism 
\[
\cO_X(U)\otimes_A\cU_A(L)
\cong\cU_{\cO_X(U)}\bigl(\cO_X(U)\otimes_A L\bigr)
\cong\cU_{\cO_X(U)}\bigl(\cT_X(U)\bigr)
\]
of \(K\)-algebras.
\end{Remark}

Let \(X\) be a rigid analytic \(K\)-space. In \cite{MNM1991} and \cite{Raczka2024Holonomic}, the sheaf \(\cD_X\) is defined as a sheafified Grothendieck construction for rings of differential operators. For an affinoid \(K\)-space \(X=\Sp A\) such that \(\cT_X\) is free, this constructions canonically agrees with the one of \autoref{Definition:Algebraic_Differential_Operators_D}. This shows that our approach is compatible with the one of \cite{Mebkhout1989} and \cite{Raczka2024Holonomic} for smooth rigid analytic \(K\)-spaces. 

In \autoref{Appendix:D_Modules_Rigid_Spaces:Coherent} we collect the basics of the theory of coherent and holonomic \(\cD\)-modules on rigid analytic \(K\)-spaces. We denote the respective categories by \(\Coh(\cD_X)\) and \(\Hol(\cD_X)\).

\subsection{Overview of analytic \texorpdfstring{\(\Dcap\)}{D-cap}-modules}\label{Subsection:Overview_Dcap} Our main concern is the study of finite length properties of the analytic \(\Dcap\)-modules developed in \cite{AW2019,AW2018,ABW2021}. We briefly summarise the aspects of these papers we need. Although there is subsequent work simplifying and refining aspects of the theory, for the purpose of this overview we mostly follow the aforementioned papers. The starting point is the construction of a sheaf of analytic differential operators \(\Dcap\), which is obtained from the sheaf of algebraic differential operators \(\cD\) by taking certain completions.

Let \(A\) be an affinoid \(K\)-algebra and let \(L\) be a smooth \((K,A)\)-Lie algebra. Let \(\cA\) be an affine formal model of \(A\) together with an \((R,\cA)\)-Lie lattice \(\cL\) of \(L\). Recall from \cite[Definition~6.2]{AW2019} that the \bemph{Fréchet completed enveloping algebra of \(\bm L\)} is
\[
\wideparen{\cU_A(L)}
\coloneqq\varprojlim\nolimits_n\widehat{\cU_{\cA}\bigl(\pi^n\cL\bigr)_K}\,.
\]
As \(\cU_{\cA}(\cL)\) is a deformable \(R\)-algebra with \(n\)-th deformation \(\cU_{\cA}(\pi^n\cL)\), the \(K\)-algebras \(\widehat{\cU_{\cA}\bigl(\pi^n\cL\bigr)_K}\) are the associated complete doubly filtered \(K\)-algebras, as discussed in \autoref{Subsection:Stabilisation_Doubly_Filtered_Algebras}. It is shown in \cite[Section~6.2]{AW2019} that \(\wideparen{\cU_A(L)}\) does not depend on the choice of the \((R,\cA)\)-Lie lattice \(\cL\). 

\begin{Definition}[{{\cite[Definition~9.3]{AW2019}}}]\label{Definition:Analytic_Differential_Operators_Dcap}
Let \(X\) be a smooth rigid analytic \(K\)-space. The sheaf \(\bm{\Dcap_X}\) is the sheaf of \(K\)-algebras defined by
\[
\Dcap_X(U)
\coloneqq\wideparen{\cU_{\cO_X(U)}\bigl(\cT_X(U)\bigr)}
\]
for \(U\) an affinoid subdomain of \(X\).
\end{Definition}

\begin{Remark}
This is a sheaf by \cite[Theorem~1.1]{Bode2019}.
\end{Remark}

The Fréchet completions appearing in \autoref{Definition:Analytic_Differential_Operators_Dcap} are typically not noetherian, although their limit constituents are. In order to obtain a well-behaved theory of \(\Dcap\)-modules, one uses the framework of Fréchet--Stein algebras and their coadmissible modules, as developed in \cite{ST2003}.

\begin{Definition}[{{\cite[Section~3]{ST2003}}}]
A \(K\)-algebra \(A\) is a \bemph{Fréchet--Stein \(\bm K\)-algebra} if \(A\cong\varprojlim_n A_n\) is a countable projective limit of noetherian Banach \(K\)-algebras \(A_n\) whose connecting morphisms have dense image and turn \(A_n\) into a (two-sided) flat \(A_{n+1}\)-module. 

An \(A\)-module \(M\) is \bemph{coadmissible} if \(M\cong\varprojlim_n M_n\) is a countable projective limit of finitely generated \(A_n\)-modules \(M_n\) such that the natural morphism \(A_n\otimes_{A_{n+1}}M_{n+1}\to M_n\) is an isomorphism.
\end{Definition}

Let \(A\) be a Fréchet--Stein \(K\)-algebra. The category of coadmissible \(A\)-modules \(\bm{\cC_A}\) is a full abelian subcategory of the category of \(A\)-modules by \cite[Corollary~3.5]{ST2003}. The subcategory \(\cC_A\) should be thought of as a replacement for the full subcategory of finitely generated \(A\)-modules. A coadmissible \(A\)-module is of \bemph{finite length} if it is an object of finite length in the abelian category \(\cC_A\).

Let \(X\) be a smooth rigid analytic \(K\)-space. It was shown in \cite[Theorem~6.3]{AW2013} that \(\Dcap_X(U)\) is a Fréchet--Stein \(K\)-algebra for affinoid subdomains \(U\) for which the \((K,\cO_X(U))\)-Lie algebra \(\cT_X(U)\) admits a smooth Lie lattice. A \(\Dcap_X\)-module \(\cM\) is \bemph{coadmissible}, if there is an admissible open cover by affinoid subdomains \(U\) for which the \((K,\cO_X(U))\)-Lie algebra \(\cT_X(U)\) admits a smooth Lie lattice and such that \(\restrict{\cM}{U}\) is associated to a coadmissible \(\cD_X(U)\)-module via a completed localisation procedure (see \cite[Section~8.2]{AW2019}). By the above and \cite[Proposition~8.3]{AW2019}, the category of coadmissible \(\Dcap_X\)-modules \(\bm{\cC_X}\) is abelian. A coadmissible \(\Dcap_X\)-module is of \bemph{finite length} if it is an object of finite length in the abelian category \(\cC_X\).

Let \(A\cong\varprojlim_n A_n\) be a Fréchet--Stein \(K\)-algebra. If there is \(d\) such that each \(A_n\) is Auslander regular of global dimension \(d_n\le d\), \cite[Section~8]{ST2003} introduces and studies a dimension function on \(\cC_A\). The condition that \(A_n\) be Auslander regular of global dimension \(d_n\le d\) was subsequently relaxed to \(A_n\) being Auslander--Gorenstein of self-injective dimension \(d_n\le d\) in \cite[Section~5]{ABW2021}. In this incarnation, this dimension theory can be applied to the theory of coadmissible \(\Dcap\)-modules.

Let \(A\) be an affinoid \(K\)-algebra and let \(L\) be a smooth \((K,A)\)-Lie algebra. Let \(\cA\) be an affine formal model of \(A\) together with an \((R,\cA)\)-Lie lattice \(\cL\) of \(L\) of rank \(r\). By \cite[Definition~5.1]{ABW2021}, the \bemph{dimension} of a (non-zero) coadmissible \(\wideparen{\cU_A(L)}\)-module \(M\) is defined by \(\rd(M)=\dim A+r-j_{\wideparen{\cU_A(L)}}(M)\). Let \(X\) be a smooth rigid analytic \(K\)-space and \(\cM\) a coadmissible \(\Dcap_X\)-module. For \(\sU\) an admissible cover by affinoid subspaces \(U\) for which \(\cT_X(U)\) admits a smooth Lie lattice, \cite[Definition~5.3]{ABW2021} defines the \bemph{dimension} of \(\cM\) by \(\rd_\sU(\cM)=\sup\{\rd(\cM(U))\,|\,U\in\sU\}\). This is independent of the chosen cover by \cite[Lemma~5.3]{ABW2021}. Moreover, it is shown in \cite[Theorem~6.2]{ABW2021} that Bernstein's inequality holds. A coadmissible \(\Dcap_X\)-module \(\cM\) is called \bemph{weakly holonomic} if \(\cM=0\) or \(\rd(\cM)=\dim X\). The category of weakly holonomic \(\Dcap_X\)-modules \(\bm{\cC_X^\wh}\) is an abelian subcategory of \(\cC_X\) by \cite[Proposition~7.1]{ABW2021}. In particular, a \(\Dcap_X\)-module is of finite length as weakly holonomic \(\Dcap_X\)-module if it is a weakly holonomic \(\Dcap_X\)-module of finite length.

\subsection{\texorpdfstring{\(\Dcap\)}{D-cap}-modules on the closed unit polydisc}\label{Subsection:Dcap_Unit_Polydisc} We now specialise to the closed unit polydisc in order to relate the theory surveyed in \autoref{Subsection:Overview_Dcap} to the completed Weyl algebras we considered in \autoref{Section:Completed_Weyl_Algebras}. Let \(A=K\langle x_1,\dots,x_d\rangle\) and let \(X=\Sp A\) be the closed unit polydisc over \(K\). The tangent sheaf \(\cT_X\) is free with global sections
\[
L
\coloneqq\cT_X(X)
=\Der_K(A)
=\bigoplus\nolimits_{i=1}^d A\,\del_i
\]
for \(\del_i\coloneqq\del_{x_i}\) the partial derivative with respect to \(x_i\) on \(A\). This is a smooth \((K,A)\)-Lie algebra. An affine formal model for \(A\) is \(\cA=R\langle x_1,\dots,x_d\rangle\) and
\[
\cL
\coloneqq\Der_R(\cA)
=\bigoplus\nolimits_{i=1}^d \cA\,\del_i
\]
defines a canonical smooth \((R,\cA)\)-Lie lattice in \(L\) of rank \(d\). In particular, \(L\) admits a smooth Lie lattice and we observe that
\[
\Dcap_X(X)
=\wideparen{\cU_A(L)}
=\varprojlim\nolimits_n\widehat{\cU_\cA(\pi^n\cL)_K}
\]
is a Fréchet--Stein presentation of \(\Dcap_X(X)\). Recall that sending a coadmissible \(\Dcap_X\)-module \(\cM\) to the coadmissible \(\Dcap_X(X)\)-module \(\cM(X)\) defines an equivalence of categories by \cite[Theorem~9.5]{AW2019}. We now express this Fréchet--Stein presentation for \(\Dcap_X(X)\) in terms of completed Weyl algebras.

\begin{Proposition}\label{Proposition:Frechet_Stein_Completed_Weyl_Algebras}
Let \(A=K\langle x_1,\dots,x_d\rangle\) and let \(\cA=R\langle x_1,\dots,x_d\rangle\). Consider the \((K,A)\)-Lie algebra \(\Der_K(A)\) and the \((R,\cA)\)-Lie lattice \(\Der_R(\cA)\). For all \(n\ge0\), there is a natural isomorphism
\[
\widehat{A_{d,n}(R)}
\longto
\widehat{\cU_{\cA}(\pi^n\cL)}
\]
compatible with the natural transition maps. In particular, there is a natural isomorphism
\[
\wideparen{A_d(K)}
\coloneqq\varprojlim\nolimits_n\widehat{A_{d,n}(K)}
\longto\varprojlim\nolimits_n\widehat{\cU_{\cA}(\pi^n\cL)_K}
=\wideparen{U_A(L)}
\]
of Fréchet--Stein \(K\)-algebras.
\end{Proposition}

\begin{Proof}
Recall from \cite[Chapter~2, p.27]{GW1989} that for a ring \(S\) and a derivation \(\delta\colon S\to S\), the skew polynomial ring \(S[y;\delta]\) is freely generated by \(y\) over \(S\) and subject to the relation \([y,s]=ys-sy=\delta(s)\) for all \(s\in S\). If \(S\) is an \(R\)-algebra and \(\delta\) and \(R\)-derivation, then \(S[y;\delta]\) is an \(R\)-algebra.

Note that the \(A\)-module \(L\) as well as the \(\cA\)-module \(\cL\) are freely generated by \(\del_1,\dots,\del_d\). From this we conclude that there is a natural isomorphism \(\cU_\cA\bigl(\pi^n\cL\bigr)\cong\cA[y_1;\pi^n\del_1]\cdots[y_d;\pi^n\del_d]\) of \(R\)-algebras for all \(n\ge0\). On the other hand, in the course of the proof of \autoref{Lemma:Completed_Weyl_Algebras_Slices} we have already observed that there is natural isomorphism \(A_{d,n}(R)\cong R[x_1,\dots,x_d][y_1;\pi^n\del_1]\cdots[y_d;\pi^n\del_d]\) of \(R\)-algebras for all \(n\ge0\).

Thus, the natural embedding \(R[x_1,\dots,x_d]\to R\langle x_1,\dots,x_d\rangle\) induces natural \(R\)-algebra morphism
\begin{equation}\label{Equation:Natural_Morphism}
A_{d,n}(R)\longto\cU_\cA\bigl(\pi^n\cL\bigr)\tag{\(\ast\)}
\end{equation}
for \(n\ge0\), which induces natural \(R/(\pi^m)\)-algebra morphisms
\[
A_{d,n}(R)/(\pi^m)\longto\cU_\cA\bigl(\pi^n\cL\bigr)/(\pi^m)
\]
for all \(m>0\). As \(R\langle x_1,\dots,x_d\rangle\) is the \(\pi\)-adic completion of \(R[x_1,\dots,x_d]\), these morphisms are isomorphisms for all \(m>0\) and all \(n\ge0\). In particular, the natural morphism \equationref{Equation:Natural_Morphism} induces an isomorphism of \(\pi\)-adic completions for all \(n\ge0\). Finally, observe that the transition maps are always induced by \(\pi^{n+1}\del_i=\pi(\pi^n\del_i)\), hence are compatible with the isomorphisms constructed.
\end{Proof}

We have the following compatibility between the dimension theory for coadmissible \(\wideparen{A_d(K)}\)-modules introduced in \autoref{Subsection:Overview_Dcap} and the dimension theories studied in \autoref{Subsection:Dimensions_Completed_Weyl_Algebras}.

\begin{Corollary}\label{Corollary:Dimensions_Frechet_Stein}
If \(\wideparen{M}\) is a coadmissible \(\wideparen{A_d(K)}\)-module and \(M_n=\widehat{A_{d,n}(K)}\otimes_{\wideparen{A_d(K)}}\wideparen{M}\), then
\[
\rd(\wideparen{M})
=\sup\bigl\{\rd_n(M_n)\,\big|\,n\ge0\bigr\}\,.
\]
\end{Corollary}

\begin{Proof}
By definition, we have \(\rd(\wideparen{M})=2d-j_{\wideparen{A_d(K)}}(\wideparen{M})\). On the other hand, for all \(n\ge0\) we know that \(\rd_n(M_n)=2d-j_{\widehat{A_{d,n}(K)}}(M_n)\) by  \autoref{Proposition:Dimensions_Equal}. For each \(n\ge0\), the \(K\)-algebra \(\widehat{A_{d,n}(K)}\) is Auslander regular of global dimension at most \(2d\) by \autoref{Corollary:Noetherian_Auslander_Regular}. Thus, as already observed in \cite[Section~8]{ST2003}, we find that \(j_{\wideparen{A_d(K)}}(\wideparen{M})=\min\bigl\{j_{\widehat{A_{d,n}(K)}}(M_n)\,\big|\,n\ge0\bigr\}\). We then conclude that \(\rd(\wideparen{M})=\sup\bigl\{\rd_n(M_n)\,\big|\,n\ge0\bigr\}\).
\end{Proof}

We need the following elementary observations.

\begin{Lemma}\label{Lemma:Coadmissible_Finite_Length}
Let \(A=\varprojlim_n A_n\) be a Fréchet--Stein \(K\)-algebra and let \(M=\varprojlim_n M_n\) be a non-zero coadmissible \(A\)-module. 
\begin{enumroman}
\item\label{Lemma:Coadmissible_Finite_Length:Non_Zero} There is an integer \(n_0\ge0\) such that \(M_n\) is non-zero for all \(n\ge n_0\).
\item\label{Lemma:Coadmissible_Finite_Length:Length} If there are integers \(n_0\ge0\) and \(\ell\ge0\) such that \(M_n\) is an \(A_n\)-module of finite length bounded by \(\ell\) for all \(n\ge n_0\), then \(M\) is of finite length bounded by \(\ell\).
\end{enumroman}
\end{Lemma}

\begin{Proof}
For \ref{Lemma:Coadmissible_Finite_Length:Non_Zero}, note that if \(M_n=0\) for all \(n\ge0\), then \(M=\varprojlim_n M_n=0\). Thus, there is \(n_0\ge0\) such that \(M_{n_0}\) is non-zero. Note that
\[
A_{n_0}\otimes_{A_{n_0+1}}A_{n_0+1}\otimes_{A_{n_0+2}}\cdots\otimes_{A_n}M_n
\cong M_{n_0}
\]
for all \(n\ge n_0\). This forces \(M_n\) to be non-zero for all \(n\ge n_0\).

For \ref{Lemma:Coadmissible_Finite_Length:Length}, consider a chain
\[
N_0\leq N_1\leq N_2\leq\cdots\leq N_r\leq M
\]
of coadmissible submodules of \(M\) such that \(N_i/N_{i-1}\) is non-zero for all \(1\le i\le r\). As there are finitely many such quotients, \ref{Lemma:Coadmissible_Finite_Length:Non_Zero} shows that there is \(n_0\ge0\) such that \(N_{i,n}/N_{i,n-1}\) is a non-zero \(A_n\)-module for all \(n\ge n_0\). Thus, for all \(n\ge n_0\) there is a chain
\[
N_{0,n}\leq N_{1,n}\leq N_{2,n}\leq\cdots\leq N_{r,n}\leq M_n
\]
of submodules of \(M_n\) such that \(N_{i,n}/N_{i,n-1}\) is non-zero for all \(1\le i\le r\). We may assume that \(M_n\) is of finite length bounded by \(\ell\) for all \(n\ge n_0\) as well. But this forces \(r\le\ell\), which shows that \(M\) is of finite length bounded by \(\ell\).
\end{Proof}

The following now provides a criterion of a coadmissible \(\wideparen{A_d(K)}\)-module to be of finite length, in terms of stabilisation conditions for the multiplicities of the associated finitely generated \(\widehat{A_{d,n}(K)}\)-modules.

\begin{Proposition}\label{Proposition:Finite_Length_Criterion}
Let \(\wideparen{M}\) be a non-zero coadmissible \(\wideparen{A_d(K)}\)-module and write \(M_n=\widehat{A_{d,n}(K)}\otimes_{\wideparen{A_d(K)}}\wideparen{M}\). Suppose that \(\wideparen{M}\) is weakly holonomic and that there are integers \(n_0\ge0\) and \(m>0\) such that \(\rm_n(M_n)\le m\) for all \(n\ge n_0\). Then \(\wideparen{M}\) is of finite length bounded by \(m\).
\end{Proposition}

\begin{Proof}
As \(\wideparen{M}\) is a non-zero weakly holonomic \(\wideparen{A_d(K)}\)-module, we observe that \(\rd(\wideparen{M})=d\). By \autoref{Corollary:Dimensions_Frechet_Stein} this shows that \(\rd_n(M_n)\le d\) for all \(n\ge0\). By \autorefpart{Lemma:Coadmissible_Finite_Length}{Lemma:Coadmissible_Finite_Length:Non_Zero}, we know that \(M_n\) is non-zero for sufficiently large \(n\). Without loss of generality, assume that \(n_0\ge0\) is such that \(M_n\) is non-zero and \(\rm_n(M_n)\le m\) for all \(n\ge n_0\). By Bernstein's inequality \autoref{Corollary:Bernstein_Inequality_Bernstein} we then conclude that \(\rd_n(M_n)=d\) for \(n\ge n_0\). Hence, for all \(n\ge n_0\), we see that \(M_n\) is a non-zero holonomic \(\widehat{A_{d,n}(K)}\)-module of multiplicity \(\rm_n(M_n)\le m\). Thus, \(M_n\) is of finite length bounded by \(m\) for all \(n\ge n_0\) by \autoref{Corollary:Multiplicity_Length_Bound} and we conclude by \autorefpart{Lemma:Coadmissible_Finite_Length}{Lemma:Coadmissible_Finite_Length:Length}.
\end{Proof}

\section{The Extension Functor and Finite Length}\label{Section:Extension}

In this section we study the extension functor from coherent \(\cD\)-modules to coadmissible \(\Dcap\)-modules, previously considered in \cite[Section]{BB2021} and \cite[Section~7.2]{ABW2021}. We explicitly check its compatibility with side-changing and pushforward along a closed embedding. We then use these compatibilities together with the results of \autoref{Section:Stabilisation} and \autoref{Section:Completed_Weyl_Algebras} to show that, for a quasi-compact smooth rigid analytic \(K\)-space, the extension functor sends holonomic \(\cD\)-modules to weakly holonomic \(\Dcap\)-modules of finite length.

\subsection{The extensions functors} Let \(X\) be a smooth rigid analytic \(K\)-space. We will use the extension functor, which produces coadmissible \(\Dcap_X\)-modules from coherent \(\cD_X\)-modules.

\begin{Definition}[{{\cite[Section~2.5]{BB2021}, \cite[Section~7.2]{ABW2021}}}] 
Let \(X\) be a smooth rigid analytic \(K\)-space. The \bemph{extension functors} are defined by
\[
E_X\colon\Coh(\cD_X)\to\cC_X
,\qquad
\cM\mapsto\Dcap_X\otimes_{\cD_X}\cM
\]
and
\[
E_X^\rr\colon\Coh(\cD_X^\op)\to\cC_X^\rr
,\qquad
\cM\mapsto\cM\otimes_{\cD_X}\Dcap_X
\]
\end{Definition}

The extension functor is exact by \cite[Lemma~4.14]{Bode2019} and faithful by \cite[Theorem~3.1]{ABW2021}. If \(\cM\) is a holonomic \(\cD_X\)-module, then \(E_X(\cM)\) is a weakly holonomic \(\Dcap_X\)-module by \cite[Proposition~7.2]{ABW2021}. We require certain compatibility results between the extensions functors and \(\cD\)-module, respectively, \(\Dcap\)-module constructions.

Recall that for a smooth rigid analytic \(K\)-space \(X\), we denote by \(\omega_X=\det\Omega_{X/K}\) the \bemph{canonical sheaf}. For both coherent \(\cD_X\)-modules and coadmissible \(\Dcap_X\)-modules, the functors \(\omega_X\otimes_{\cO_X}\blank\) and \(\iHom_{\cO_X}(\omega_X,\blank)\) defined a canonical equivalence, called \bemph{side-changing equivalence}, between left and right modules, cf.\ \cite[Proposition~1.2.12]{HTT2008} and \cite[Theorem~3.5]{AW2018}.

For \(i\colon Y\to X\) the inclusion of a smooth closed analytic subset of a smooth rigid analytic \(K\)-space \(X\), for both coherent \(\cD\)-modules and coadmissible \(\Dcap\)-modules, there is a pushforward functor \(i_+\), cf.\ \cite[Section~2.5]{Raczka2024Holonomic} and \cite[Section~4.5]{AW2018}. For coadmissible \(\Dcap\)-modules, the pushforward functor \(i_+\) induces a version of Kashiwara's equivalence \cite[Theorem~7.1]{AW2018}.

\begin{Proposition}\label{Proposition:Extension_Functor_Compatibilities}
Let \(X\) be a smooth rigid analytic \(K\)-space and let \(i\colon Y\to X\) be the inclusion of a smooth closed analytic subset.
\begin{enumroman}
\item\label{Proposition:Extension_Functor_Compatibilities:Side_Changing} For every coherent left \(\cD_X\)-module \(\cM\), there is a natural isomorphism
\[
\omega_X\otimes_{\cO_X} E_X(\cM)
\cong E_X^\rr\bigl(\omega_X\otimes_{\cO_X}\cM\bigr)
\]
of coadmissible right \(\Dcap_X\)-modules. For every coherent right \(\cD_X\)-module \(\cM\), there is a natural isomorphism
\[
\iHom_{\cO_X}\bigl(\omega_X,E_X^\rr(\cM)\bigr)
\cong E_X\bigl(\iHom_{\cO_X}(\omega_X,\cM)\bigr)
\]
of coadmissible left \(\Dcap_X\)-modules.
\item\label{Proposition:Extension_Functor_Compatibilities:Closed_Pushforward} For every coherent left \(\cD_Y\)-module \(\cM\), there is a natural isomorphism
\[
i_+\bigl(E_Y(\cM)\bigr)
\cong E_X\bigl(i_+(\cM)\bigr)
\]
of coadmissible left \(\Dcap_X\)-modules.
\end{enumroman}
\end{Proposition}

\begin{Proof}
Note that the extension functor is compatible with restrictions. Therefore we may assume that \(\cM\) is finitely presented as a \(\cD_X\)-module. Note that \(E_X\) is exact by \cite[Lemma~4.14]{Bode2019} and for both coherent \(\cD\)-modules and coadmissible \(\Dcap\)-modules side-changing is exact as it induces an equivalence. Thus, it suffices to check \ref{Proposition:Extension_Functor_Compatibilities:Side_Changing} for the coherent left \(\cD_X\)-module \(\cD_X\).

To see \ref{Proposition:Extension_Functor_Compatibilities:Side_Changing}, we observe that
\[
\omega_X\otimes_{\cO_X}E_X\bigl(\cD_X\bigr)
=\omega_X\otimes_{\cO_X}\bigl(\Dcap_X\otimes_{\cD_X}\cD_X\bigr)
\cong\omega_X\otimes_{\cO_X}\Dcap_X
\]
while
\[
E_X^\rr(\omega_X\otimes_{\cO_X}\cD_X)
=\bigl(\omega_X\otimes_{\cO_X}\cD_X\bigr)\otimes_{\cD_X}\Dcap_X\,.
\]
Note that in the latter case, \(\omega_X\otimes_{\cO_X}\cD_X\) carries the right \(\cD_X\)-module structure as obtained by tensoring the right \(\cD_X\)-module \(\omega_X\) with the left \(\cD_X\)-module \(\cD_X\) (cf.\ \cite[Proposition~1.2.9]{HTT2008}). Up to \(\cO_X\)-linear isomorphism, this right \(\cD_X\)-module is induced by the right \(\cD_X\)-module structure on \(\cD_X\) by \cite[Lemma~3.2]{AW2018}. We conclude by the associativity of the tensor product. Since the pair \(\omega_X\otimes_{\cO_X}\blank\) and \(\iHom_{\cO_X}(\omega_X,\blank)\) induces the side-changing equivalence, the other claim of \ref{Proposition:Extension_Functor_Compatibilities:Side_Changing} follows.

By \ref{Proposition:Extension_Functor_Compatibilities:Side_Changing} it suffices to check \ref{Proposition:Extension_Functor_Compatibilities:Closed_Pushforward} for coherent \emph{right} \(\cD_Y\)-modules. Note that for coherent \(\cD_Y\)-modules \(i_+\) is right-exact, as it obtained from tensoring with the transfer bimodule \(\cD_{Y\to X}\) and then applying \(i_\ast\), which is exact by \cite[Theorem~A.1]{AW2018}. Thus, it suffices to check \ref{Proposition:Extension_Functor_Compatibilities:Closed_Pushforward} for the coherent right \(\cD_Y\)-module \(\cD_Y\). Note that then
\[
i_+\bigl(\cD_Y\bigr)
\cong i_\ast\bigl(\cD_{Y\to X}\bigr)
=i_\ast\bigl(i^\ast\cD_X\bigr)
\cong\cD_X/\cI_Y\cD_X
\]
for \(\cI_Y\) the ideal sheaf of \(Y\) in \(X\). In particular, we find that
\[
E_X^\rr\biggl(i_+\bigl(\cD_Y\bigr)\biggr)
\cong E_X^\rr\bigl(\cD_X/\cI_Y\cD_X\bigr)
\cong\Dcap_X/\cI_Y\Dcap_X\,.
\]
By \cite[Theorem~6.2]{AW2018} we can reduce to the case where
\[
i_+\bigl(\Dcap_Y\bigr)(X)
=\Dcap_Y(Y)\wideparen{\otimes}_{\Dcap_Y(Y)}\Dcap_X(X)/\cI_Y(X)\Dcap_X(X)
\cong\Dcap_X(X)/\cI_Y(X)\Dcap_X(X)
\]
as coadmissible \(\Dcap_Y(Y)\)-modules. Since, \((\Dcap_X/\cI_Y\Dcap_X)(X)\cong\Dcap_X(X)/\cI_Y(X)\Dcap_X(X)\), we are done by \cite[Theorem~9.5]{AW2019}.
\end{Proof}

\subsection{Finite length for extended \texorpdfstring{\(\Dcap\)}{D-cap}-modules}

Let \(X\) be a quasi-compact smooth rigid analytic \(K\)-space. We now show that \(E_X\) sends holonomic \(\cD_X\)-modules to weakly holonomic \(\Dcap_X\)-modules of finite length. We first establish a restricted version of this statement for the closed unit polydisc, using the stabilisation results of \autoref{Section:Stabilisation} and our considerations in \autoref{Section:Completed_Weyl_Algebras}. We then generalise to general smooth affinoid \(K\)-spaces using Kashiwara's equivalence \cite[Theorem~7.1]{AW2018}. Finally, we obtain the desired version by a simple covering argument.

\begin{Theorem}\label{Theorem:Finite_Length_Polydisc}
Let \(X=\Sp A\) with \(A=K\langle x_1,\dots,x_d\rangle\) be the closed unit polydisc and let \(\cM\) be a holonomic globally finitely presented \(\cD_X\)-module. Then \(E_X(\cM)\) is of finite length in \(\cC_X^\wh\).
\end{Theorem}

\begin{Proof}
We may assume that \(\cM\) is non-zero. Then \(E_X(\cM)\) is weakly holonomic by \cite[Proposition~7.2]{ABW2021} and non-zero by \cite[Theorem~3.1]{ABW2021}. Moreover, \(E_X(\cM)\) is of finite length if and only if \(E_X(\cM)(X)\) is of finite length by \cite[Theorem~9.5]{AW2019}. Since \(\cM\) is globally finitely presented, we conclude that \(E_X(\cM)(X)\cong\Dcap_X(X)\otimes_{\cD_X(X)}\cM(X)\). Note that \(\cM(X)\) is a finitely generated \(\cD_X(X)\)-module by \autoref{Proposition:Weak_Kiehl}.

Write \(M\coloneqq\cM(X)\) and \(\wideparen{M}\coloneqq E_X(\cM)(X)\). For the \((K,A)\)-Lie algebra \(L=\Der_K(L)\), choose the affine formal model \(\cA=R\langle x_1,\dots,x_d\rangle\) and the smooth \((R,\cA)\)-Lie lattice \(\cL=\Der_R(\cA)\). Write \(U=\cU_\cA(\cL)\), \(U_K=\cU_A(L)\) and \(\wideparen{U}=\wideparen{\cU_L(A)}\). Let \(m_1,\dots,m_r\in M\) be a finite generating set and let \(N\) be the \(\pi\)-torsionfree \(U\)-submodule of \(M\) generated by \(m_1,\dots,m_r\). Note that \(\wideparen{M}\cong\wideparen{U}\otimes_{U_K} M\cong\wideparen{U}\otimes_UN\). 

As \(U\) is a deformable \(R\)-algebra, \autoref{Proposition:Stabilisation_Doubly_Filtered_Algebra} shows that there exists an integer \(n_0\ge0\) so that
\[
\Gr^{\cB_\bullet}\bigl(\widehat{U_{n_0,K}}\otimes_UN\bigr)
\cong\Gr^{\cB_\bullet}\bigl(\widehat{U_{n,K}}\otimes_UN\bigr)
\]
for all \(n\ge n_0\). In \autoref{Proposition:Frechet_Stein_Completed_Weyl_Algebras} we saw that
\(\widehat{U_{n,K}}\cong\widehat{A_{d,n}(K)}\) and \(\wideparen{U}\cong\wideparen{A_d(K)}\). Hence, we conclude that
\[
\Gr^{\cB_\bullet}\bigl(\widehat{A_{d,n_0}(K)}\otimes_{\wideparen{A_d(K)}}\wideparen{M}\bigr)
\cong\Gr^{\cB_\bullet}\bigl(\widehat{A_{d,n_0}(K)}\otimes_UN\bigr)
\cong\Gr^{\cB_\bullet}\bigl(\widehat{A_{d,n}(K)}\otimes_UN\bigr)
\cong\Gr^{\cB_\bullet}\bigl(\widehat{A_{d,n}(K)}\otimes_{\wideparen{A_d(K)}}\wideparen{M}\bigr)
\]
for all \(n\ge n_0\). Recall that this isomorphism identifies the \(\Gr^{\cB_\bullet}\bigl(\widehat{U_{n_0,K}}\otimes_UN\bigr)\)-module structure on \(\Gr^{\cB_\bullet}\bigl(\widehat{A_{d,n_0}(K)}\otimes_{\wideparen{A_d(K)}}\wideparen{M}\bigr)\) with the \(\Gr^{\cB_\bullet}\bigl(\widehat{U_{n,K}}\otimes_UN\bigr)\)-module structure on \(\Gr^{\cB_\bullet}\bigl(\widehat{A_{d,n}(K)}\otimes_{\wideparen{A_d(K)}}\wideparen{M}\bigr)\). In particular, the characteristic ideals of \(\Gr^{\cB_\bullet}\bigl(\widehat{A_{d,n_0}(K)}\otimes_{\wideparen{A_d(K)}}\wideparen{M}\bigr)\) and \(\Gr^{\cB_\bullet}\bigl(\widehat{A_{d,n}(K)}\otimes_{\wideparen{A_d(K)}}\wideparen{M}\bigr)\) are isomorphic for all \(n\ge n_0\). By \autoref{Definition:Hilbert_Polynomial_Completed_Weyl_Algebras}, we see that
\[
\rm_{n_0}\bigl(\widehat{A_{d,n_0}(K)}\otimes_{\wideparen{A_d(K)}}\wideparen{M}\bigr)
=\rm_n\bigl(\widehat{A_{d,n}(K)}\otimes_{\wideparen{A_d(K)}}\wideparen{M}\bigr)
\]
for all \(n\ge n_0\). Set \(m\coloneqq \rm_{n_0}\bigl(\widehat{A_{d,n_0}(K)}\otimes_{\wideparen{A_d(K)}}\wideparen{M}\bigr)\). As \(\wideparen{M}\) is non-zero and weakly holonomic, we conclude by \autoref{Proposition:Finite_Length_Criterion} since \(\rm_n\bigl(\widehat{A_{d,n}(K)}\otimes_{\wideparen{A_d(K)}}\wideparen{M}\bigr)=m\) for all \(n\ge n_0\).
\end{Proof}

We record the following elementary lemma, which we need to pass from local finite length statements to global finite length statements.

\begin{Lemma}\label{Lemma:Finite_Length_Cover}
Let \(X\) be a quasi-compact smooth rigid analytic \(K\)-space. Let \(\cM\) be a weakly holonomic \(\Dcap_X\)-module. If there is an admissible open cover \(\sU=\{U_\alpha\}\) by affinoid subdomains \(U_\alpha\) such that \(\restrict{\cM}{U_\alpha}\) is of finite length in \(\cC_{U_\alpha}^\wh\), then \(\cM\) is of finite length in \(\cC_X^\wh\).
\end{Lemma}

\begin{Proof}
Since \(X\) is quasi-compact, we may assume that the given admissible open cover is given by finitely many affinoid subdomains \(U_1,\dots,U_r\). Suppose that \(\restrict{\cM}{U_i}\) is of finite length \(\ell_i\). Let \(\ell=\max\{\ell_1,\dots,\ell_r\}\). Then each \(\restrict{\cM}{U_i}\) is of finite length at most \(\ell\). We claim that \(\cM\) is of finite length at most \(\ell r\). If not, there is a strictly ascending chain
\[
\cM_0\lneq\cM_1\lneq\cdots\lneq\cM_{\ell r+1}
\]
of submodules. In particular, each \(\cM_j/\cM_{j-1}\) is non-zero for \(1\le j\le \ell r+1\). As the \(U_1,\dots,U_r\) form an admissible open cover, this means that there is at least one \(i\) so that \(\restrict{\bigl(\cM_j/\cM_{j-1}\bigr)}{U_i}\cong\restrict{\cM_j}{U_i}/\restrict{\cM_{j-1}}{U_i}\) is non-zero. Since there are \(\ell r+1\) quotients of this form, but only \(r\) affinoid subdomains in the admissible open cover, there is at least one \(U_i\) for which at least \(\ell+1\) quotients are non-zero. This contradicts that the \(\restrict{\cM}{U_i}\) are of finite length uniformly bounded by \(\ell\).
\end{Proof}

\begin{Theorem}\label{Theorem:Finite_Length_Affinoid}
Let \(X=\Sp A\) be a smooth affinoid \(K\)-space and let \(\cM\) be a holonomic \(\cD_X\)-module. Then \(E_X(\cM)\) is of finite length in \(\cC_X^\wh\).
\end{Theorem}

\begin{Proof}
By \cite[Proposition~7.2]{ABW2021} we know that \(E_X(\cM)\) is a weakly holonomic \(\Dcap_X\)-module. Using Kashiwara's equivalence \cite[Theorem~7.1]{AW2018} we reduce to \autoref{Theorem:Finite_Length_Polydisc}.

First, assume that \(X=\Sp A\) admits a coordinate system and that \(\cM\) is a holonomic globally finitely presented \(\cD_X\)-module. Fix a closed embedding \(i\colon X\to Y\) with \(Y=\Sp K\langle x_1,\dots,x_d\rangle\) a closed unit polydisc. By \autorefpart{Proposition:Extension_Functor_Compatibilities}{Proposition:Extension_Functor_Compatibilities:Closed_Pushforward}, we know that \(E_Y\bigl(i_+(\cM)\bigr)\cong i_+\bigl(E_X(\cM)\bigr)\). Since \(X\) and \(Y\) admit coordinate systems and \(\cM\) is globally finitely presented, we conclude that \(i_+(\cM)\) is globally finitely presented by \autoref{Lemma:Pushforward_Gloablly_Finitely_Presented}. Moreover, \(i_+(\cM)\) is holonomic by \autorefpart{Lemma:Pushforward_Coherent_Holonomic}{Lemma:Pushforward_Coherent_Holonomic:Holonomic}. Thus, \(E_Y\bigl(i_+(\cM)\bigr)\) is of finite length in \(\cC_Y^\wh\) by \autoref{Theorem:Finite_Length_Polydisc}. As Kashiwara's equivalence \cite[Theorem~7.1]{AW2018} shows that \(i_+\) is a fully faithful exact embedding, we conclude that \(E_X(\cM)\) is of finite length in \(\cC_X^\wh\).

Now let \(X=\Sp A\) be an arbitrary smooth affinoid \(K\)-space and let \(\cM\) be a holonomic \(\cD_X\)-module. By definition, there is an admissible open cover \(\sU=\{U_\alpha\}\) by smooth affinoid subdomains \(U_\alpha\) admitting coordinate systems such that \(\restrict{\cM}{U_\alpha}\) is globally finitely presented on \(U_\alpha\) and \(\cM(U_\alpha)\) is a holonomic \(\cD_{U_\alpha}(U_\alpha)\)-module. Equivalently, by \autoref{Proposition:Globally_Finitely_Presented_Holonomic}, \(\restrict{\cM}{U_\alpha}\) is holonomic. By the above, note that \(E_{U_\alpha}\bigl(\restrict{\cM}{U_\alpha}\bigr)\) is of finite length in \(\cC_{U_i}^\wh\). As \(\restrict{E_X(\cM)}{U_\alpha}\cong E_{U_\alpha}\bigl(\restrict{\cM}{U_\alpha}\bigr)\) we conclude by \autoref{Lemma:Finite_Length_Cover}.
\end{Proof}

\begin{Remark}
In \cite[Theorem~V.1.3]{Milicic} a similar strategy is used to show that holonomic \(\cD\)-modules on a smooth complex variety are of finite length by reducing the case of affine varieties to the case of affine space via Kashiwara's equivalence.
\end{Remark}

\begin{Corollary}\label{Corollary:Finite_Length_Quasi_Compact}
Let \(X\) be a quasi-compact smooth rigid analytic \(K\)-space and let \(\cM\) be a holonomic \(\cD_X\)-module. Then \(E_X(\cM)\) is of finite length in \(\cC_X^\wh\).
\end{Corollary}

\begin{Proof}
By \cite[Proposition~7.2]{ABW2021} we know that \(E_X(\cM)\) is a weakly holonomic \(\Dcap_X\)-module. Consider an admissible open cover \(\sU=\{U_\alpha\}\) by affinoid subdomains \(U_\alpha\) so that each \(\restrict{\cM}{U_\alpha}\) is a holonomic \(\cD_{U_\alpha}\)-module. Thus, \(E_{U_\alpha}\bigl(\restrict{\cM}{U_\alpha}\bigr)\) is weakly holonomic \(\Dcap_{U_\alpha}\)-module of finite length by \autoref{Theorem:Finite_Length_Affinoid}. Now we conclude by \autoref{Lemma:Finite_Length_Cover}.
\end{Proof}

\begin{Remark}
After the results of this paper were obtained, we were informed by Raoul Hallopeau that he obtained the same result for quasi-compact smooth formal curves over \(R\), see \cite[Exemple~6.4]{Hallopeau2025Formal}. Together with \cite[Proposition~6.1]{Hallopeau2025Rigid} this provides an alternative proof of \autoref{Corollary:Finite_Length_Quasi_Compact} for quasi-compact smooth rigid analytic curves over \(K\) with good reduction.
\end{Remark}

\section{Geometric Examples of \texorpdfstring{\(\Dcap\)}{Dcap}-Modules of Finite Length}\label{Section:Finite_Length}

Using our previous examination of the extension functor, we strengthen results of \cite{BB2021} and \cite{ABW2021} on the coadmissibility of \(\Dcap\)-modules. These proofs follow a uniform pattern: for a suitable \(\Dcap_X\)-module \(\cM\) of interest, one associates a coherent \(\cD_X\)-module \(\cN\) which is equipped with a natural morphism \(E_X(\cN)\to\cM\) of \(\Dcap_X\)-modules with dense image. Coadmissibility of \(\cM\) is equivalent to this morphism being a surjection, which can then be used to infer more refined information about \(\cM\) by studying \(\cN\). We use this principle to deduce various finite length results from \autoref{Theorem:Finite_Length_Affinoid} and \autoref{Corollary:Finite_Length_Quasi_Compact}.

\subsection{Examples obtained from meromorphic connections} Let \(X=\Sp A\) be a smooth affinoid \(K\)-space admitting a coordinate system. Let \(f\in A\) be non-constant and consider the closed analytic subset \(Z=\{f=0\}\). The \bemph{sheaf of meromorphic functions with poles in \(\bm Z\)}, denote \(\bm{\cO_X(\ast Z)}\), is defined by \(\cO_X(\ast Z)(U)=B[f^{-1}]\) for \(U=\Sp B\) an affinoid subdomain of \(X\). Note that \(\cO_X(\ast Z)\) is naturally a \(\cD_X\)-module. Consider the sheaf of rings \(\bm{\cD_X(\ast Z)}\coloneqq\cO_X(\ast Z)\otimes_{\cO_X}\cD_X\). A \bemph{meromorphic connection with singularities along \(\bm Z\)} is an \(\cD_X(\ast Z)\)-module which is coherent as \(\cO_X(\ast Z)\)-module.

As in the case of coherent \(\cD_X\)-modules, there are some subtleties when working with coherent \(\cO_X(\ast Z)\)-modules (cf.\ our remark between \autoref{Definition:Globally_Finitely_Presented_Coherent} and \autoref{Proposition:Weak_Kiehl}). In particular, it is not known if the global sections of a coherent \(\cO_X(\ast Z)\)-module are finitely generated over \(\cO_X(\ast Z)(X)=A[f^{-1}]\). However, restricting to meromorphic connections with singularities along \(Z\) which are globally finitely presented as \(\cO_X(\ast Z)\)-modules resolves these issues. In this restricted form, we can strengthen \cite[Theorem~4.4]{BB2021} (note that these subtleties are not considered loc.\ cit.).

\begin{Theorem}\label{Theorem:BB_Extension}
Let \(X=\Sp A\) be a smooth affinoid \(K\)-space admitting a coordinate system. Let \(f\in A\) be non-constant. Define \(Z=\{f=0\}\) and let \(U=X\setminus Z\). Denote the open embedding by \(j\colon U\to X\). Let \(\cN\) be a meromorphic connection on \(X\) with singularities along \(Z\) which is globally finitely presented as \(\cO_X(\ast Z)\)-module and let \(\cM=E_U\bigl(\restrict{\cN}{U}\bigr)\). Let \(m_1,\dots,m_k\) be a finite generating set of \(\cN(X)\) as \(A[f^{-1}]\)-module, and let \(b_1,\dots,b_k\) be the corresponding \(b\)-functions from \cite[Théorème~3.1.1]{MNM1991}.

If all roots of the \(b_i\) are of positive type (in the sense of \cite[Definition~13.1.1]{Kedlaya2022}), then \(j_\ast(\cM)\) is of finite length in \(\cC_X^\wh\).
\end{Theorem}

\begin{Proof}
Given our assumptions, \cite[Theorem~4.4]{BB2021} applies and shows that the natural morphism \(E_X(\cN)\to j_\ast(\cM)\) is an isomorphism of coadmissible \(\Dcap_X\)-modules. In particular, \(j_\ast(\cM)\) is weakly holonomic by \cite[Proposition~7.2 and Proposition~7.1]{ABW2021}. It now suffices to show that \(\cN\) is a holonomic \(\cD_X\)-module to conclude using \autoref{Theorem:Finite_Length_Affinoid}. Since we assume \(\cN\) to be globally finitely presented as \(\cO_X(\ast Z)\)-module, we may restrict to working with its global sections \(\cN(X)\). Note that \(\cN(X)\) is a holonomic \(\cD_X(X)\)-module by \cite[Théorème~3.1.1]{MNM1991}, hence \(\cN\) is a holonomic \(\cD_X\)-module by \autoref{Proposition:Globally_Finitely_Presented_Holonomic}. This concludes the proof.
\end{Proof}

\begin{Corollary}\label{Corollary:Non_Holonomic_Finite_Length}
Let \(X=K\langle x\rangle\) be the \(1\)-dimensional closed unit disc. Let \(Z=\{x=0\}\) and \(U=X\setminus Z\). Write \(j\colon U\to X\). Then there exists a non-holonomic weakly holonomic \(\Dcap_X\)-module of finite length.
\end{Corollary}

\begin{Proof}
In \cite[Lemma~6.2]{Bode2025Holonomic}, it is shown that \(\cM\coloneqq\Dcap_X/\Dcap_X(x\del-\lambda)\) for \(\lambda\in K\) of positive type such that \(p\lambda\) is of type \(0\) is a non-holonomic weakly holonomic \(\Dcap_X\)-module satisfying \(\cM\cong j_\ast\restrict{\cM}{U}\). Since \(\lambda\) is assumed to be of positive type, \autoref{Theorem:BB_Extension} shows that \(j_\ast\bigl(\restrict{\cM}{U}\bigr)\) is of finite length.
\end{Proof}

\subsection{Examples obtained from integrable connections} Recall that an \bemph{integrable connection} is an \(\cO_X\)-coherent \(\cD_X\)-module. By \cite[Theorem~7.3]{AW2018}, this is equivalent to being an \(\cO_X\)-coherent coadmissible \(\Dcap_X\)-module. In \cite{ABW2021}, coadmissibility and weak holonomicity of \(\Dcap_X\)-modules obtained from integrable connections is studied. We first obtain the following strengthening of \cite[Lemma~10.5]{ABW2021}.

\begin{Lemma}\label{Lemma:ABW_Distinguished_Open_Extension}
Let \(X=\Sp A\) be a smooth affinoid \(K\)-space and let \(f\in A\) be non-constant. Define \(Z=\{f=0\}\) and let \(U=X\setminus Z\). Denote the embedding by \(j\colon U\to X\). If \(\cM\) is an integrable connection on \(X\), then \(j_\ast\bigl(\restrict{\cM}{U}\bigr)\) is a of finite length in \(\cC_X^\wh\).
\end{Lemma}

\begin{Proof}
We know that \(j_\ast\bigl(\restrict{\cM}{U}\bigr)\) is weakly holonomic by \cite[Lemma~10.5]{ABW2021}. To conclude that \(j_\ast\bigl(\restrict{\cM}{U}\bigr)\) is of finite length, we proceed as in the proof \cite[Lemma~10.5]{ABW2021}. By \autoref{Lemma:Finite_Length_Cover} we may assume that \(X\) admits a coordinate system. In particular, by \cite[Theorem~9.5]{AW2019} it suffices to show that \(\cM(U)=j_\ast\bigl(\restrict{\cM}{U}\bigr)(X)\) is of finite length as \(\Dcap_X(X)\)-module.  Denote \(M=\cM(X)\). By \cite[Théorème~3.2.1]{MNM1991}, \(M[f^{-1}]\) is a holonomic \(\cD_X(X)\)-module. Since we know that \(j_\ast\bigl(\restrict{\cM}{U}\bigr)\) is coadmissible, the natural morphism
\[
\Dcap_X(X)\otimes_{\cD_X(X)} M[f^{-1}]\longto\cM(U)
\]
is surjective by \cite[Proposition~2.14]{BB2021}. Note that the domain of the surjection can be understood as \(E_X(\cM(\ast Z))\). As \(\cM(\ast Z)\) is a holonomic \(\cD_X\)-module by \cite[Théorème~4.3.3]{MNM1991}, we conclude that \(E_X(\cM(\ast Z))\) is weakly holonomic \(\Dcap_X\)-module of finite length by \autoref{Theorem:Finite_Length_Affinoid}. In particular, \(\Dcap_X(X)\otimes_{\cD_X(X)} M[f^{-1}]\) is of finite length as \(\Dcap_X(X)\)-module and hence so is \(\cM(U)\).
\end{Proof}

Based on \autoref{Lemma:ABW_Distinguished_Open_Extension}, we now obtain the following strengthening of \cite[Theorem~10.5]{ABW2021}.

\begin{Theorem}\label{Theorem:ABW_Open_Extension}
Let \(U\) be a Zariski-open subset of a quasi-compact smooth rigid analytic \(K\)-space. Denote the open embedding by \(j\colon U\to X\). If \(\cM\) is an integrable connection on \(X\), then \(\rR^q j_\ast\bigl(\restrict{\cM}{U}\bigr)\) is of finite length in \(\cC_X^\wh\) for all \(q\ge0\).
\end{Theorem}

\begin{Proof}
As \(X\) is assumed to be quasi-compact, we can reduce to the affinoid case using \autoref{Lemma:Finite_Length_Cover}. We now proceed as in the proof of \cite[Theorem~10.5]{ABW2021}. Let \(X=\Sp A\) be a smooth affinoid \(K\)-space and \(U=X\setminus\cV(f_1,\dots,f_r)\) for \(f_1,\dots,f_r\in A\) non-constant. The admissible open subsets \(X_{f_i}=X\setminus\cV(f_i)\) form an admissible open cover \(\cU\) and we have \(\rR^q j_\ast\bigl(\restrict{\cM}{U}\bigr)\cong\check{\rH}^q\bigl(\cU,\restrict{\cM}{U}\bigr)\) for all \(q\ge0\) by \cite[Corollary~10.1]{ABW2021}. The \(q\)-th cochain group of the \v{C}ech complex has the form
\[
\prod_{i_0,\dots,i_q}\restrict{\cM}{U}\bigl(X_{f_{i_0}}\cap\dots\cap X_{f_{i_q}}\bigr)
=\prod_{i_0,\dots,i_q}\cM\bigl(X_{f_{i_0}\cdots f_{i_q}}\bigr)
\]
and each \(\cM\bigl(X_{f_{i_0}\cdots f_{i_q}}\bigr)\) is \(\cD_X\)-module of finite length by \autoref{Lemma:ABW_Distinguished_Open_Extension}. As the cover is finite, each \(q\)-cochain group is a finite product of \(\cD_X\)-modules of finite length, hence of finite length as well. In particular, \(\check{\rH}^q\bigl(\cU,\restrict{\cM}{U}\bigr)\) is a \(\Dcap_X\)-module of finite length for all \(q\ge0\). This concludes the proof.
\end{Proof}

Finally, we obtain a strengthening of \cite[Theorem~10.6]{ABW2021}. Recall that for \(Z\) a closed analytic subset of a smooth rigid analytic \(K\)-space \(X\) and \(\cF\) a coherent \(\cO_X\)-module, we can define for all \(q\ge0\) a sheaf \(\underline\cH_Z^q(\cF)\) called the \bemph{local cohomology sheaf of \(\bm\cF\) with support in \(\bm Z\)}. Here \(\underline\cH_Z^q(\blank)\) is the \(q\)-th derived functor of \(\underline\cH_Z^0(\blank)\), which assigns to a coherent \(\cO_X\)-module \(\cF\) the maximal subsheaf with support in \(Z\) (cf.\ \cite[Definition~2.1.3]{Kisin1999}).

\begin{Theorem}\label{Theorem:ABW_Closed_Extension}
Let \(Z\) be a closed analytic subset of a quasi-compact smooth rigid analytic \(K\)-space. If \(\cM\) is an integrable connection on \(X\), then \(\underline\cH_Z^q(\cM)\) is of finite length in \(\cC_X^\wh\) for all \(q\ge0\).
\end{Theorem}

\begin{Proof}
As in the proof \autoref{Theorem:ABW_Open_Extension} we may reduce to the affinoid case using \autoref{Lemma:Finite_Length_Cover}. Now we proceed as in the proof of \cite[Theorem~10.6]{ABW2021}. Denote by \(j\colon U\to X\) the embedding of the Zariski-open complement of \(Z\). As argued in the proof of \cite[Theorem~10.6]{ABW2021}, there are isomorphisms \(\underline\cH_Z^q(\cM)\cong\rR^{q-1}j_\ast\bigl(\restrict{\cM}{U}\bigr)\) for all \(q\ge2\) and an exact sequence
\[
\begin{tikzcd}
0 \ar[r] & \underline\cH_Z^0(\cM) \ar[r] & \cM \ar[r] & j_\ast\bigl(\restrict{\cM}{U}\bigr) \ar[r] & \underline\cH_Z^1(\cM) \ar[r] & 0\,.    
\end{tikzcd}
\]
Note \(\cM\) is a \(\Dcap_X\)-module of finite length since \(\cM\) is an \(\cO_X\)-coherent coadmissible \(\Dcap_X\)-module by \cite[Theorem~7.3]{AW2018}. Moreover, \(\rR^q j_\ast\bigl(\restrict{\cM}{U}\bigr)\) is a \(\Dcap_X\)-module of finite length by \autoref{Theorem:ABW_Open_Extension} for all \(q\ge0\). This concludes the proof.
\end{Proof}

\appendix

\section{Recollection on Hilbert Polynomials}\label{Appendix:Hilbert_Polynomials}

Fix a base field \(F\). Note that we neither assume that \(F\) is algebraically closed nor that \(F\) is of characteristic \(0\). We recall the construction of the Hilbert polynomials associated to finitely generated modules over almost commutative \(F\)-algebras, following \cite{MCR2001}.

We recall some terminology from \cite[Section~8.1.9]{MCR2001}. An \(F\)-algebra \(A\) is \bemph{affine} if its finitely generated as an \(F\)-algebra by a finite-dimensional \(F\)-subspace \(V\). Let \(A\) be a filtered \(F\)-algebra with filtration \(\cF_\bullet A\) and let \(M\) be a filtered \(A\)-module with filtration \(\cF_\bullet M\). The filtration \(\cF_\bullet A\) is \bemph{standard} if \(\cF_i A=(\cF_1)^i\) for each \(i\) and \bemph{finite-dimensional} if \(\cF_0 A=F\) and \(\dim_F\cF_iA<\infty\) for each \(i\). The filtration \(\cF_\bullet M\) is \bemph{standard} if \(\cF_i M=\cF_i A\cdot\cF_0 M\) for each \(i\) and \bemph{finite-dimensional} if \(\dim_F\cF_iM<\infty\) for each \(i\). The filtration \(\cF_\bullet M\) is \bemph{good} if \(\gr^{\cF_\bullet} M\) is a finitely generated \(\gr^{\cF_\bullet} A\)-module.

\begin{Definition}[{{\cite[Section~8.4.2]{MCR2001}}}]\label{Definition:Somewhat_Commutative}
Let \(A\) be an \(F\)-algebra. We say that \(A\) is \bemph{almost commutative} if it an affine \(F\)-algebra which has some standard finite-dimensional filtration \(\cF_\bullet A\) such that \(\gr^{\cF_\bullet}A\) is a commutative.
\end{Definition}

If \(A\) is an almost commutative \(F\)-algebra via the filtration \(\cF_\bullet A\) and \(M\) is finitely generated \(A\)-module with a good filtration \(\cF_\bullet M\), then \(\cF_\bullet M\) is finite-dimensional since \(\cF_\bullet A\) is (standard) finite-dimensional.

\begin{Definition}\label{Definition:Hilbert_Function}
Let \(A\) be an almost commutative \(F\)-algebra via the filtration \(\cF_\bullet A\) and let \(M\) be a finitely generated \(A\)-module with a good filtration \(\cF_\bullet M\). The function
\[
h_{M,\cF_\bullet}\colon\bN_0\longto\bN_0
,\,\qquad
i\mapsto\dim_F\cF_iM
\]
is called the \bemph{Hilbert function of \(\bm{(M,\cF_{\bullet} M)}\)}.
\end{Definition}

The following is the central result.

\begin{Theorem}[{{\cite[Theorem~8.4.5,~Section~8.6.11]{MCR2001}}}]\label{Theorem:Hilbert_Polynomial}
Let \(A\) be an almost commutative \(F\)-algebra via the filtration \(\cF_\bullet A\) and let \(M\) be a finitely generated \(A\)-module with good filtration \(\cF_\bullet M\). Then there is a polynomial \(p_{M,\cF_\bullet}\in\bQ[x]\) such that
\[
h_{M,\cF_\bullet}(i)
=p_{M,\cF_\bullet}(i)
\]
for all \(i\gg0\).
\end{Theorem}

\begin{Remark}
In \cite[Section~8.4]{MCR2001} filtrations on finitely generated \(A\)-modules, for \(A\) almost commutative, are fixed to be standard finite-dimensional. Only in \cite[Section~8.6]{MCR2001}, dealing with the more general notion of a \emph{somewhat commutative \(F\)-algebra} (see \cite[Paragraph~8.6.9]{MCR2001}), the generality of good filtrations is discussed. As any almost commutative \(F\)-algebra is somewhat commutative, the results of \cite[Section~8.6]{MCR2001} apply and we will refer to them as appropriate.
\end{Remark}

We recall the following proposition from \cite{MCR2001}.

\begin{Proposition}[{{\cite[Proposition~8.4.4]{MCR2001}}}]\label{Proposition:Numerical_Functions}
Let \(f\colon\bN\to\bQ\) be a function.
\begin{enumalph}
\item The following conditions, for \(i\gg0\), are equivalent:
\begin{enumroman}\label{Proposition:Itemizea}
\item\label{Proposition:Numerical_Functions:Polynomial} \(f(i)=p(i)\) for some \(p\in\bQ[x]\) with \(\deg p=d\).
\item\label{Proposition:Numerical_Functions:Polynomial_Form} \(f(i)=a_d\binom id+\cdots+a_1\binom i1+a_0\) with \(a_i\in\bQ\) and \(a_d\ne0\).
\item \(f(i+1)-f(i)=a_d\binom i{d-1}+\cdots+a_2\binom i1+a_1\) with \(a_i\in\bQ\) and \(a_d\ne0\).
\end{enumroman}
\item If the conditions in \ref{Proposition:Itemizea} are satisfied, then the following hold:
\begin{enumroman}\label{Proposition:Itemizeb}
\item \(p\) and \(a_0,\dots,a_d\) are uniquely determined by \(f\);
\item if \(f(i)\in\bZ\) for \(i\gg0\), then each \(a_0,\dots,a_d\in\bZ\);
\item\label{Proposition:Numerical_Functions:Nonnegative_Coefficients} if \(f(i)\in\bN_0\) for \(i\gg0\), then each \(a_d\in\bN_0\).
\end{enumroman}
\end{enumalph}
\end{Proposition}

By \autorefparts{Proposition:Numerical_Functions}{Proposition:Itemizea}{Proposition:Numerical_Functions:Polynomial}, the polynomial \(p_{M,\cF_\bullet}\) of \autoref{Theorem:Hilbert_Polynomial} is uniquely determined by \(h_{M,\cF_\bullet}\).

\begin{Definition}\label{Definition:Hilbert_Polynomial}
The polynomial \(p_{M,\cF_\bullet}\in\bQ[x]\) from \autoref{Theorem:Hilbert_Polynomial} is called a \bemph{Hilbert polynomial of \(\bm M\)}.
\end{Definition}

Note that the Hilbert function \(h_{M,\cF_\bullet}\) and the Hilbert polynomial \(p_{M,\cF_\bullet}\) satisfy \autorefparts{Proposition:Numerical_Functions}{Proposition:Itemizea}{Proposition:Numerical_Functions:Polynomial}. Thus, by \autorefparts{Proposition:Numerical_Functions}{Proposition:Itemizea}{Proposition:Numerical_Functions:Polynomial_Form} we can write
\[
p_{M,\cF_\bullet}(x)
=\sum_{j=0}^d a_j\binom xj
\,,\qquad\binom xj\coloneqq\frac{x(x-1)\cdots(x-j+1)}{j!}\,,
\]
as rational polynomials are determined by their values for \(i\gg0\). Now \autorefparts{Proposition:Numerical_Functions}{Proposition:Itemizeb}{Proposition:Numerical_Functions:Nonnegative_Coefficients} shows that \(a_d\in\bN_0\). 

By \cite[Theorem~8.6.19]{MCR2001}, if we fix the filtration \(\cF_\bullet A\) of an almost commutative \(F\)-algebra \(A\), then for any finitely generated \(A\)-module \(M\) the degree and the leading term of a Hilbert polynomial of \(M\) are independent of the choice of a good filtration \(\cF_\bullet M\). This prompts the following definition.

\begin{Definition}\label{Definition:Dimension_Multiplicity}
Let \(A\) be an almost commutative \(F\)-algebra via the filtration \(\cF_\bullet A\) and let \(M\) be a finitely generated \(A\)-module. Let \(\cF_\bullet M\) be good filtration with associated Hilbert polynomial \(h_{M,\cF_\bullet}(x)=\sum_{j=0}^d a_j\binom xj\). The \bemph{dimension of \(\bm M\)} is \(\bm{\rd(M)}\coloneqq\deg h_{M,\cF_\bullet}=d\) and the \bemph{multiplicity of \(\bm M\)} is \(\bm{\rm(M)}\coloneqq a_d\).
\end{Definition}

\begin{Proposition}[{{\cite[Theorem~8.4.8, Corollary~8.6.20]{MCR2001}}}]\label{Proposition:Dimension_Multiplicity_SES}
Let \(A\) be an almost commutative \(F\)-algebra via the filtration \(\cF_\bullet A\) and let 
\[
\begin{tikzcd}
0 \ar[r] & M' \ar[r] & M \ar[r] & M'' \ar[r] & 0
\end{tikzcd}
\]
be a a short exact sequence of finitely generated \(A\)-modules.   
\begin{enumroman}
\item\label{Proposition:Dimension_Multiplicity_SES:Dimension} \(\rd(M)=\max\big(\rd(M'),\rd(M'')\big)\).
\item\label{Proposition:Dimension_Multiplicity_SES:Cases} Precisely one of the following statements is true
\begin{enumalph}
\item \(\rd(M')<\rd(M)=\rd(M'')\) and \(\rm(M)=\rm(M'')\).
\item \(\rd(M'')<\rd(M)=\rd(M')\) and \(\rm(M)=\rm(M')\).
\item \(\rd(M')=\rd(M)=\rd(M'')\) and \(\rm(M)=\rm(M')+\rm(M'')\).    
\end{enumalph}
\end{enumroman}
\end{Proposition}

In our situation of interest, the dimension introduced in \autoref{Definition:Dimension_Multiplicity} can be understood geometrically. For the lack of an adequate reference for the level of generality we require, we provide a short proof. 

\begin{Proposition}\label{Proposition:Hilbert_Polynomial_Degree_Dimension}
Let \(A\) be an almost commutative \(F\)-algebra via the filtration \(\cF_\bullet A\). Let \(M\) be a finitely generated \(A\)-module and fix a good filtration \(\cF_\bullet M\). Then
\[
\dim V\big(\Ann \gr^{\cF_\bullet} M\big)
=\deg h_{M,\cF_\bullet}
\]
where \(V\big(\Ann \gr^{\cF_\bullet} M\big)\subseteq\Spec\gr^{\cF_\bullet} A\).
\end{Proposition}

\begin{Proof}
Let \(B=\bigoplus_{i\ge0}B_i\) be a graded finitely generated \(F\)-algebra generated by \(B_1\) and let \(N=\bigoplus_{i\ge0} N_i\) be a finitely generated graded \(S\)-module. As discussed in \cite[Example~23.93]{GW2023}, we can consider \(X=\Proj B\) and the coherent \(\cO_X\)-module \(\cN=\widetilde{N}\) associated to \(N\). Then \cite[Proposition and Definition~23.91]{GW2023} applies and shows that there is a polynomial \(\Phi_\cN\in\bQ[x]\) such that
\[
\dim_F N_i
=\Phi_\cN(i)
\]
for all \(i\gg0\) and \(\deg\Phi_\cN=\dim\Supp\cN\). Since \(\cN\) is the coherent \(\cO_X\)-module associated to the graded \(B\)-module \(N\), we observe that \(\Supp\cN=V_+(\Ann N)\) as subsets of \(X=\Proj B\). This makes sense as \(\Ann N\) is homogeneous since \(N\) is graded.

Let us now specialise to \(N=\gr^{\cF_\bullet}M\) with \(N_i=\gr_i^{\cF_\bullet}M\). Then \(\dim_F\gr_i^{\cF_\bullet}M=h_{M,\cF_\bullet}(i)-h_{M,\cF_\bullet}(i-1)\) for all \(i\). Moreover, \(\Supp\cN=V_+(\Ann\gr^{\cF_\bullet}M)\) for \(\cN=\widetilde{\gr^{\cF_\bullet}M}\). From the equivalences in \autorefpart{Proposition:Numerical_Functions}{Proposition:Itemizea} and the above we conclude that \(\deg h_{M,\cF_\bullet}=\dim V_+(\Ann\gr^{\cF_\bullet}M)+1\). On the other hand, since \(V(\Ann\gr^{\cF_\bullet}M)\) is the affine cone of \(V_+(\Ann\gr^{\cF_\bullet}M)\), it follows that \(V(\Ann\gr^{\cF_\bullet}M)=\dim V_+(\Ann\gr^{\cF_\bullet}M)+1\). Thus,
\[
\dim V\big(\Ann \gr^{\cF_\bullet} M\big)
=\dim V_+(\Ann\gr^{\cF_\bullet}M)+1
=\deg h_{M,\cF_\bullet}
\]
as claimed.
\end{Proof}

\section{Algebraic \texorpdfstring{\(\cD\)}{D}-Modules on Rigid Analytic Spaces}\label{Appendix:D_Modules_Rigid_Spaces}

\subsection{Smoothness and coordinate systems}\label{Appendix:D_Modules_Rigid_Spaces:Coordinates} Let \(X\) be a rigid analytic \(K\)-space. Recall from \cite[Example~4.4.1]{FvdP2004} the \bemph{sheaf of differentials \(\bm{\Omega_X}\)}, which is determined by \(\Omega_X(U)=\Omega_{\cO_X(U)/K}^f\) for any affinoid subdomain \(U\) of \(X\). Here \(\Omega_{\cO_X(U)/K}^f\) denotes the \bemph{universal finite differential module of \(\bm{\cO_X(U)}\)}, as defined in \cite[Section~3.6]{FvdP2004}. The dual of \(\Omega_X\) is called the \bemph{tangent sheaf} and denoted by \(\bm{\cT_X}\). We say that \(X\) is \bemph{smooth} if \(\Omega_X\) is locally free.

To argue locally with \(\cD\)-modules, the affinoid subdomains considered need to be particularly well-behaved. We compare different notions for an affinoid \(K\)-space \(X=\Sp A\). In \cite[Section~2]{Raczka2024Holonomic}, \(X\) is said to \bemph{admit a (global) coordinate system}, if there exist \(f_1,\dots,f_d\in A\) such that \(\rd f_1,\dots,\rd f_d\in\Omega_{A/K}\) is an \(A\)-basis. In \cite[Section~4.3.1]{MNM1991}, \(X\) is called \bemph{carte affino\"ide} if it admits a coordinate system. In \cite[Section~6.3]{Bode2026}, \(X\) is said to \bemph{admit a local coordinate system} if there exist \(f_1,\dots,f_d\in A\) and an \(A\)-basis \(\del_1,\dots,\del_d\in\Der_K(A)\) such that \([\del_i,\del_j]=0\) and \(\del_i(f_j)=\delta_{ij}\). It is immediate (cf.\ \cite[Section~2.3.2]{Raczka2024Holonomic}) that a coordinate system defines a local coordinate system. In fact, both notions are equivalent to admitting an étale morphism (cf.\ \cite[Definition~2.1]{BLR1995}) to a polydisc of dimension \(d\).

\begin{Lemma}\label{Lemma:Coordinate_Systems_Equivalence}
Let \(X=\Sp A\) be an affinoid \(K\)-space. The following are equivalent:
\begin{enumroman}
\item\label{Lemma:Coordinate_Systems_Equivalence:Étale} There exists an étale morphism \(X\to\Sp K\langle x_1,\dots,x_d\rangle\);
\item\label{Lemma:Coordinate_Systems_Equivalence:Raczka} \(X\) admits a coordinate system;
\item\label{Lemma:Coordinate_Systems_Equivalence:Bode} \(X\) admits a local coordinate system;
\end{enumroman}
\end{Lemma}

\begin{Proof}
We show \ref{Lemma:Coordinate_Systems_Equivalence:Étale}\(\iff\)\ref{Lemma:Coordinate_Systems_Equivalence:Raczka} and \ref{Lemma:Coordinate_Systems_Equivalence:Raczka}\(\iff\)\ref{Lemma:Coordinate_Systems_Equivalence:Bode}.

A morphism of \(K\)-spaces \(\varphi\colon X\to\Sp K\langle x_1,\dots,x_d\rangle\) corresponds to a morphism \(\varphi^\sharp\colon K\langle x_1,\dots,x_d\rangle\to A\) of the underlying affinoid \(K\)-algebras, which is completely determined by \(f_i=\varphi^\sharp(x_i)\). If \(\varphi\) is étale, the natural morphism \(A\otimes_{K\langle x_1,\dots,x_d\rangle}\Omega_{K\langle x_1\dots,x_d\rangle/K}\to\Omega_{A/K}\)
is an isomorphism of \(A\)-modules, mapping \(\rd x_i\) to \(\rd f_i\). In particular, the elements \(\rd f_1,\dots,\rd f_r\) form an \(A\)-basis for \(\Omega_{A/K}\). Conversely, any \(f_1,\dots,f_d\in A\) such that \(\rd f_1,\dots,\rd f_d\in\Omega_{A/K}\) form an \(A\)-basis defines a natural morphism for which sending \(\rd x_i\) to \(\rd f_i\) defines an isomorphism \(A\otimes_{K\langle x_1,\dots,x_d\rangle}\Omega_{K\langle x_1\dots,x_d\rangle/K}\to\Omega_{A/K}\) of \(A\)-modules.

For \ref{Lemma:Coordinate_Systems_Equivalence:Raczka}\(\implies\)\ref{Lemma:Coordinate_Systems_Equivalence:Bode}, denote by \(\del_1,\dots,\del_d\in\Der_K(A)\) the dual basis to \(\rd f_1,\dots,\rd f_d\). By definition, we then have \(\del_i(f_j)=\del_i(\rd f_j)=\delta_{ij}\). Moreover, writing \([\del_i,\del_j]=\sum_{i=1}^d g_k^{ij}\del_k\) and evaluating at the \(f_1,\dots,f_d\) shows that \([\del_i,\del_j]=0\). For \ref{Lemma:Coordinate_Systems_Equivalence:Bode}\(\implies\)\ref{Lemma:Coordinate_Systems_Equivalence:Raczka}, consider the elements \(\rd f_1,\dots,\rd f_d\in\Omega_{A/K}\). As \(\del_i(\rd f_j)=\del_i(f_j)=\delta_{ij}\), we conclude that \(\rd f_1,\dots,\rd f_d\in\Omega_{A/K}\) is the dual basis to \(\del_1,\dots,\del_d\).
\end{Proof}

\begin{Proposition}
Let \(X=\Sp A\) be a smooth affinoid \(K\)-space. Then \(X\) admits a coordinate system if and only if \(\Omega_X\) is free if and only if \(\cT_X\) is free.
\end{Proposition}

\begin{Proof}
If \(\Omega_X\) is free, then \(\cT_X\) is free, being its dual. Now assume that \(\cT_X\) is free. On global sections, this shows that \(\Der_K(A)\) is a free \(A\)-module and it remains to show that we can find a local coordinate system. Arguing locally at maximal ideals, this is the content of \cite[Theorem~99]{Matsumura1980} (cf.\ \cite[Section~1.1.2]{MNM1991} and \cite[Lemma~2.3]{Raczka2024Holonomic}). Finally, suppose that \(X\) admits a coordinate system. By \autoref{Lemma:Coordinate_Systems_Equivalence} there exits an étale morphism \(\varphi\colon X\to\Sp K\langle x_1,\dots,x_d\rangle\). In particular, there is an isomorphism \(\varphi^\ast\Omega_{\Sp K\langle x_1,\dots,x_d\rangle}\cong\Omega_X\) by \cite[Proposition~2.6]{BLR1995}, which shows that \(\Omega_X\) is free. 
\end{Proof}

Let \(X=\Sp A\) be an affinoid \(K\)-space. Following \cite[Definition~9.3]{AW2013}, we say that the \((K,A)\)-Lie algebra \(\Der_K(A)\) \bemph{admits a smooth Lie lattice}, if there is an affine formal model \(\cA\) for \(A\) and a smooth \((R,\cA)\)-Lie lattice \(\cL\) in \(\Der_K(A)\). The following is now immediate from the proof of \cite[Lemma~9.3]{AW2019}.

\begin{Corollary}
Let \(X=\Sp A\) be a smooth affinoid \(K\)-space admitting a coordinate system. Then \(\Der_K(A)\) admits a smooth Lie lattice.
\end{Corollary}

\subsection{Coherent and holonomic \texorpdfstring{\(\cD\)}{D}-modules}\label{Appendix:D_Modules_Rigid_Spaces:Coherent} Let \(X\) be a smooth rigid analytic \(K\)-space. In \autoref{Definition:Algebraic_Differential_Operators_D}, we introduced the sheaf \(\cD_X\) of algebraic differential operators. In this section, we will go over some fundamentals of the theory of coherent and holonomic \(\cD_X\)-modules.

\begin{Definition}[{{\cite[Section~4.3.1]{MNM1991}, \cite[Section~2.2]{Raczka2024Holonomic}}}]\label{Definition:Globally_Finitely_Presented_Coherent}
Let \(X\) be a smooth rigid analytic \(K\)-space and let \(\cM\) be a \(\cD_X\)-module.
\begin{enumroman}
\item We say that \(\cM\) is \bemph{globally finitely presented}, if it admits a finite presentation
\[
\begin{tikzcd}
\cD_X^s \ar[r] & \cD_X^r \ar[r] & \cM \ar[r] & 0\,.    
\end{tikzcd}
\]
\item We say that \(\cM\) is \bemph{coherent}, if there exists an admissible covering \(\sU=\{U_\alpha\}\) by open affinoid subdomains admitting coordinate systems such that for every \(U_\alpha\), the restriction \(\restrict{\cM}{U_\alpha}\) is a globally finitely presented \(\cD_{U_i}\)-module.
\end{enumroman}
\end{Definition}

Let \(X=\Sp A\) be an affinoid \(K\)-space. By a theorem of Kiehl \cite[Theorem~9.4.3/3]{BGR1984}, taking global sections induces an equivalence between coherent \(\cO_X\)-modules and finitely generated \(\cO_X(X)\)-modules. If \(X\) admits a coordinate system, then \cite[Theorem~9.5]{AW2019} provides the analogues statement for coadmissible \(\Dcap_X\)-modules and coadmissible \(\Dcap_X(X)\)-modules. However, in this generality, it is not known whether this theorem holds for coherent \(\cD_X\)-modules and finitely generated \(\cD_X(X)\)-modules (cf.\ \cite[Section~2.4]{Raczka2024Holonomic}). Consequently, handling coherent \(\cD_X\)-modules becomes a bit awkward and requires some care. However, we still have the following, which suffices for our purposes.

\begin{Proposition}[{{\cite[Lemma~2.5(3)]{Raczka2024Holonomic}}}]\label{Proposition:Weak_Kiehl}
Let \(X=\Sp A\) be a smooth affinoid \(K\)-space admitting a coordinate system. Taking global sections induces an equivalence between the category of globally finitely presented \(\cD_X\)-modules and the category of finitely generated \(\cD_X(X)\)-modules.
\end{Proposition}

For \(X=\Sp A\) a smooth affinoid \(K\)-space admitting a coordinate system, the \(K\)-algebra \(\cD_X(X)\) satisfies the assumptions of \cite[Section~1.1.2]{MNM1991} (cf.\ \cite[Section~4.3]{MNM1991} and \cite[Section~2.3]{Raczka2024Holonomic}). In particular, we can define what it means for a finitely generated \(\cD_X(X)\)-module to be \bemph{holonomic} (or of \bemph{minimal dimension}), cf.\ \cite[Definition~1.2.4]{MNM1991} and \cite[Section~2.3]{Raczka2024Holonomic}.

\begin{Definition}[{{\cite[Section~4.3.1]{MNM1991}, \cite[Section~2.3]{Raczka2024Holonomic}}}]
Let \(X\) be a smooth rigid analytic \(K\)-space and let \(\cM\) be a \(\cD_X\)-module. We say that \(\cM\) is \bemph{holonomic}, if there exists an admissible covering \(\sU=\{U_\alpha\}\) by open affinoid subdomains admitting coordinate systems such that for every \(U_\alpha\), the restriction \(\restrict{\cM}{U_\alpha}\) is a globally finitely presented \(\cD_{U_i}\)-module and the \(\cD_{U_\alpha}(U_\alpha)\)-module \(\cM(U_\alpha)\) is holonomic.
\end{Definition}

\begin{Proposition}[{{\cite[Lemma~2.7]{Raczka2024Holonomic}}}]\label{Proposition:Globally_Finitely_Presented_Holonomic}
Let \(X=\Sp A\) be a smooth affinoid \(K\)-space admitting a coordinate system and let \(\cM\) be a globally finitely presented \(\cD_X\)-module. Then \(\cM\) is holonomic if and only if \(\cM(X)\) is a holonomic \(\cD_X(X)\)-module.
\end{Proposition}

The following result is crucial.

\begin{Lemma}[{{\cite[Proposition~4.2(1)]{Raczka2024Holonomic}}}]\label{Lemma:Pushforward_Gloablly_Finitely_Presented}
Let \(i\colon Y\to X\) be the inclusion of a smooth closed analytic subset of a smooth affinoid \(K\)-space \(X\). If both \(Y\) and \(X\) admit coordinate systems and \(\cM\) is a globally finitely presented \(\cD_Y\)-module, then \(i_+(\cM)\) is a globally finitely presented \(\cD_X\)-module.
\end{Lemma}

\begin{Lemma}[{{\cite[Proposition~4.2(2),(3)]{Raczka2024Holonomic}}}]\label{Lemma:Pushforward_Coherent_Holonomic}
Let \(i\colon Y\to X\) be the inclusion of a smooth closed analytic subset of a smooth affinoid \(K\)-space \(X\). 
\begin{enumroman}
\item\label{Lemma:Pushforward_Coherent_Holonomic:Coherent} If \(\cM\) is a coherent \(\cD_Y\)-module, then \(i_+(\cM)\) is a coherent \(\cD_X\)-module.
\item\label{Lemma:Pushforward_Coherent_Holonomic:Holonomic} If \(\cM\) is a holonomic \(\cD_Y\)-module, then \(i_+(\cM)\) is a holonomic \(\cD_X\)-module.
\end{enumroman}
\end{Lemma}

\bibliographystyle{alpha}
\bibliography{references}

\end{document}